\documentclass[10pt]{article}
\usepackage{amsthm, amsmath, amsfonts,	xcolor,	hyperref,	enumerate,	geometry, multirow,	caption}
\theoremstyle{thmstyleone}
\newtheorem{theorem}{Theorem}
\newtheorem{proposition}[theorem]{Proposition} 
\newtheorem{lemma}[theorem]{Lemma}
\newtheorem{corollary}[theorem]{Corollary}
\theoremstyle{thmstyletwo}
\newtheorem{example}[theorem]{Example}
\newtheorem{remark}[theorem]{Remark}
\theoremstyle{thmstylethree}
\newtheorem{definition}[theorem]{Definition}

\title{\textbf{General Casorati inequalities and implications for Riemannian maps and Riemannian submersions}}

\author{Ravindra Singh, Kiran Meena, Kapish Chand Meena}

\date{}

\begin{document}
	\maketitle
	\begin{abstract}
		\noindent This paper presents general forms of Casorati inequalities for Riemannian maps and Riemannian submersions between Riemannian manifolds. Using these general forms, we obtain Casorati inequalities for Riemannian maps (resp. submersions) whose target (resp. source) spaces are generalized complex and generalized Sasakian space forms. As a consequence, we give Casorati inequalities for Riemannian maps (resp. submersions) when the target (resp. source) spaces are real, complex, real K\"ahler, Sasakian, Kenmotsu, cosymplectic, and almost $C(\alpha)$ space forms. To support these general forms, in the particular cases when the target or source spaces are real, complex, Sasakian, and Kenmotsu space forms, we verify known Casorati inequalities for Riemannian maps and Riemannian submersions. Further, we give Casorati inequalities for invariant and anti-invariant Riemannian maps (resp. submersions) whose target (resp. source) spaces are generalized complex and generalized Sasakian space forms. Toward information on geometric characteristics, we discuss the equality cases. We also exemplify the general forms.
	\end{abstract}
	
	\noindent \textbf{Keywords: }{Riemannian manifolds, space forms, Casorati and scalar curvatures, Riemannian submersions, Riemannian maps.}\\
	
	\noindent \textbf{MSC Classification:} {53B20, 53B35, 53C15, 53D15}
	
	\section{Introduction}
	In differential geometry, Gaussian curvature is a fundamental tool because it provides information about the shape of a surface and its points. However, it may vanish even on curved surfaces. To solve this issue, Casorati \cite{Casorati_1890} introduced the concept of \textit{Casorati curvature} for regular surfaces, which is an extrinsic invariant, useful for the visualization of shapes and appearances, and vanishes only at planar points \cite{DHV_2008, OV_2011}. Thereafter, various geometers obtained inequalities involving normalized Casorati curvatures (now known as \textit{Casorati inequalities}) for various submanifolds of ambient spaces such as real, complex, generalized complex, cosymplectic, Sasakian, Kenmotsu, generalized Sasakian, etc. (see \cite{ALVY, ASJ, DHV_2008, Ghisoiu, Lone_2017a, Vilcu, LLV_2017, LLV_2020, LLV_2022, Lone_2017, Lone_2019, Lone_2019a, Siddiqui_2018, Zhang_Zhang, Zhang_Pan_Zhang}). The Casorati inequalities were also obtained for the submanifolds of other ambient spaces, such as quasi-real, quaternionic, generalized $(\kappa, \mu)$, statistical space forms, etc. For a detailed survey, we refer to \cite{Chen_Survey}. 
	\vspace{0.2cm}
	
	Riemannian submersions are smooth maps and have applications in various fields, including physics, mechanics, relativity, spacetime, robotics, supergravity, superstring, Kaluza-Klein, Yang-Mills theories, etc. The concept of a Riemannian map generalizes that of a submanifold and a Riemannian submersion \cite{Fischer_1992}. It has rich geometry and applications \cite{GRK_book, Sahin_book}, and satisfies the eikonal equation (a bridge between geometric and physical optics). In addition, Riemannian maps allow us to compare the geometric properties of the source and target Riemannian manifolds. Due to the interesting properties mentioned above, Riemannian submersions and Riemannian maps were investigated with various structures (see \cite{Falcitelli_2004, Sahin_book, Erken_2021, MSS} and references therein).
	\vspace{0.2cm}
	
	Recently, the Casorati inequalities for Riemannian submersions and Riemannian maps were obtained with the real and complex space forms in \cite{LLSV} and with the Sasakian space forms in \cite{PLS}. Such inequalities were also obtained for Riemannian maps to Kenmotsu and nearly K\"ahler space forms in \cite{Zaidi_Shanker} and \cite{Fatima}, respectively. It was concluded that these inequalities may provide a basis for bridges between geometry and physics. Motivatingly, in this paper, we obtain general forms of Casorati inequalities for Riemannian maps and Riemannian submersions, and show that these general inequalities give relationships between normalized scalar and normalized Casorati curvatures (see Theorems \ref{General_Inq_Thm_RM}, \ref{General_Inq_Thm_VD_RS}, and \ref{General_Inq_Thm_HD_RS}). These general forms are helpful in obtaining Casorati inequalities for these smooth mappings with other space forms, similar to submanifold theory. For example, we use these general forms to obtain Casorati inequalities with generalized complex and generalized Sasakian space forms (see Theorems \ref{thm_inq_gcsc_RM}, \ref{thm_inq_gssc_RM}, \ref{gcsc_rs_vert}, \ref{gssc_rs_vert}, \ref{gcsc_rs_hor}, and \ref{gssc_rs_hor}). By direct implications, we also obtain Casorati inequalities for Riemannian maps and Riemannian submersions with real, complex, real K\"ahler, Sasakian, Kenmotsu, cosymplectic, and almost $C(\alpha)$ space forms (see Corollaries \ref{cor_rm_1}, \ref{cor_rm_2}, \ref{cor_vrs_1}, \ref{cor_vrs_2}, \ref{cor_hrs_1}, and \ref{cor_hrs_2}). In addition, in particular cases (real, complex, Sasakian and Kenmotsu space forms), we verify the obtained inequalities with the inequalities given in \cite{LLSV}, \cite{PLS} and \cite{Zaidi_Shanker} (see Remarks \ref{rem_1}, \ref{rem_3}, \ref{rem_inv_map_kenmotsu}, \ref{rem_2}, \ref{rem_inv_submersion}, and \ref{rem_4}). Further, we obtain Casorati inequalities for invariant and anti-invariant cases (see Propositions \ref{thm_inv_rm_1}, \ref{thm_inv_rm_2}, \ref{thm_inv_rs_1}, and \ref{thm_inv_rs_2}). We also discuss the equality cases and observe that they give information about some geometric characteristics.\\
	
	First, we recall the notions of the space forms mentioned above.
	
	\begin{definition} $\cite{Olszak, Tricerri, Yano}$
		Let $({M}, g, J)$ be an even-dimensional Riemannian manifold with an almost complex structure $J$, $J^{2}=-I$, the identity map. Then it is called an {\it almost Hermitian manifold} if, for all ${\cal Z}_{1}, 
		{\cal Z}_{2}\in \Gamma (T{M})$, 
		\[
		{g}(J{\cal Z}_{1}, J{\cal Z}_{2})=g({\cal Z}_{1}, {\cal Z}_{2}). 
		\]
		Furthermore, it is called a generalized	complex space form (denoted by ${M}(c_{1}, c_{2})$) if there exist smooth functions $c_{1}$ and $c_{2}$ on ${M}$ such that for all ${\cal Z}_{1}, {\cal Z}_{2}, {\cal Z}_{3}\in \Gamma (T{M})$ its curvature tensor $R$ satisfies
		\begin{align}\label{curvature_for_gcsf}
			R({\cal Z}_{1}, {\cal Z}_{2}){\cal Z}_{3}& =c_{1}\{g({\cal Z}_{2}, 			{\cal Z}_{3}){\cal Z}_{1}-g({\cal Z}_{1}, {\cal Z}_{3}){\cal Z}_{2}\} 			\nonumber \\			& +c_{2}\{g({\cal Z}_{1}, J{\cal Z}_{3})J{\cal Z}_{2}-g({\cal Z}_{2}, J			{\cal Z}_{3})J{\cal Z}_{1}+2g({\cal Z}_{1}, J{\cal Z}_{2})J{\cal Z}_{3}\}.
		\end{align}
	\end{definition}
	\noindent In particular, we also have the following notions \cite{Vanhecke, Yano, Tricerri}.
	\begin{enumerate}
		\item If $c_1=c$ and $c_2=0$, then $M(c)$ is a real space form.
		\item If $c_1=c_2=\frac{c}{4}$, then $M(c)$ is a complex space form.
		\item If $c_1=\frac{c+3\alpha}{4}$ and $c_2=\frac{c-\alpha}{4}$, then $M(c)$ is a real K\"ahler space form. 
	\end{enumerate}
	
	\begin{definition} $\cite{ABC, Blair}$
		Let $({M}, \phi, \xi, \eta, g)$ be an odd-dimensional Riemannian manifold with a $(1, 1)$ tensor field $\phi $, a structure vector field $\xi$, and $1$-form $\eta $. Then it is called an {\it almost contact metric manifold} if, for all ${\cal Z}, {\cal Z}_{1}, {\cal Z}_{2}\in \Gamma (T{M})$, 
		\[
		\eta (\xi)=1, ~\phi^{2}({\cal Z})=-{\cal Z}+\eta ({\cal Z})\xi ~{\rm and~}{g}(\phi {\cal Z}_{1}, \phi {\cal Z}_{2})=g({\cal Z}_{1}, {\cal Z}_{2})-\eta ({\cal Z}_{1})\eta ({\cal Z}_{2}). 
		\]
		In this case, we also have $\phi \xi =0$ and $\eta \circ \phi =0$. Furthermore, it is called a generalized Sasakian space form (denoted by ${M}(c_{1}, c_{2}, c_{3})$) if there exist smooth functions $c_{1}$, $c_{2}$ and $c_{3}$ on ${M}$ such that for all ${\cal Z}_{1}, {\cal Z}_{2}, {\cal Z}_{3}\in \Gamma (T{M})$ its curvature tensor $R$ satisfies
		\begin{align}\label{curvature_for_gssf}
			R({\cal Z}_{1}, {\cal Z}_{2}){\cal Z}_{3} &=c_{1}\{g({\cal Z}_{2}, {\cal Z}_{3}){\cal Z}_{1}-g({\cal Z}_{1}, 			{\cal Z}_{3}){\cal Z}_{2}\} \nonumber \\& +c_{2}\{g({\cal Z}_{1}, \phi {\cal Z}_{3})\phi {\cal Z}_{2}-g({\cal Z}_{2}, \phi {\cal Z}_{3})\phi {\cal Z}_{1}+2g({\cal Z}_{1}, \phi {\cal Z}_{2})\phi {\cal Z}_{3}\} \nonumber \\& +c_{3}\{\eta ({\cal Z}_{1})\eta ({\cal Z}_{3}){\cal Z}_{2}-\eta ({\cal Z}_{2})\eta ({\cal Z}_{3}){\cal Z}_{1}+g({\cal Z}_{1}, {\cal Z}_{3})\eta (			{\cal Z}_{2})\xi -g({\cal Z}_{2}, {\cal Z}_{3})\eta ({\cal Z}_{1})\xi\}.
		\end{align}
	\end{definition}
	
	\noindent In particular, we also have the following notions \cite{ABC, Blair}. 
	\begin{enumerate}
		\item If $c_1=\frac{c+3}{4}$ and $c_2=c_3=\frac{c-1}{4}$, then $M(c)$ is a Sasakian space form.
		\item If $c_1=\frac{c-3}{4}$ and $c_2=c_3=\frac{c+1}{4}$, then $M(c)$ is a Kenmotsu space form.
		\item If $c_1=c_2=c_3=\frac{c}{4}$, then $M(c)$ is a cosymplectic space form.
		\item If $c_1=\frac{c+3\alpha^2}{4}$ and $c_2=c_3=\frac{c-\alpha^2}{4}$, then $M(c)$ is an almost $C(\alpha)$ space form. 
	\end{enumerate}
	
	\vspace{0.2cm}
	We also use the following lemma in the further sections.
	
	\begin{lemma}\label{Lemma_Tripathi} $\cite{Tripathi_2017}$ 
		Let $\Lambda =\{(z_{1}, \dots, z_{r})\in {\Bbb R}^{r}:z_{1}+\cdots +z_{r}=k\}$ be a hyperplane of ${\Bbb R}^{r}$, and $f:{\Bbb R}^{r}\to {\Bbb R}$ be a quadratic form given by
		\begin{equation*}
			f(z_{1}, \dots, z_{r})= \lambda_1 \sum_{i=1}^{r-1}z_{i}^{2} + \lambda_2 z_{r}^{2}-2\sum_{1 \leq i<j \leq r} z_{i} z_{j}, \quad \lambda_1, \lambda_2 \in \mathbb{R}^{+}.
		\end{equation*}
		Then the constrained extremum problem $\min \limits_{(z_{1}, \dots, z_{r})\in \Lambda}f$ has the global solution $z_{1} = z_{2} = \cdots = z_{r-1} = \frac{k}{\lambda_1 + 1}, \quad z_{r} = \frac{k}{\lambda_2 + 1} = \frac{k \left( r-1\right)}{\left(\lambda_1 + 1 \right)\lambda_2} = \frac{k\left( \lambda_1 - r + 2 \right)}{\lambda_1 + 1}$, provided that $\lambda_2 = \frac{r-1}{\lambda_1 - r + 2}$.
	\end{lemma}
	
	\section{General Casorati inequalities and implications for Riemannian maps}\label{sec_2}
	This section is divided into three subsections. The first subsection is dedicated to the general form of Casorati inequalities for Riemannian maps between Riemannian manifolds. The second and third subsections are dedicated to the implications: Casorati inequalities when the target spaces are generalized complex and generalized Sasakian space forms, respectively. Now we start with some preliminaries.\\
	
	Let $F: (M_1^{m_1}, g_1) \to (M_2^{m_2}, g_2)$ be a smooth map between Riemannian manifolds with rank $3 \leq r \leq \min\{m_1, m_2\}$, and let $F_{\ast p}: T_p M_1 \to T_{F(p)} M_2$ be its derivative map at $p$. Denoting the kernel space of $F_\ast$ at $p \in M_1$ by $(\ker F_{\ast p})$ and its orthogonal complementary space in $T_p M_1$ by $(\ker F_{\ast p})^\perp$, we decompose at $p \in M_1$
	\begin{equation*}
		T_p M_1 = (\ker F_{\ast p}) \oplus (\ker F_{\ast p})^\perp.
	\end{equation*}
	Similarly we decompose at $F(p) \in M_2$,
	\begin{equation*}
		T_{F(p)} M_2 = ({\textrm{\normalfont range }} F_{\ast p}) \oplus ({\textrm{\normalfont range }} F_{\ast p})^\perp.
	\end{equation*}
	Then, the map $F$ is called a \textit{Riemannian map}, if for all $X, Y \in \Gamma(\ker F_\ast)^\perp$ the following equation holds \cite{Fischer_1992}:
	\begin{equation}\label{hor_range_isometry}
		g_1(X, Y) = g_2(F_\ast X, F_\ast Y).
	\end{equation}
	
	\noindent The bundle {\rm Hom}$\left(T{M_1}, F^{-1}T{M_2}\right)$ admits an induced connection $\nabla$ from the Levi-Civita connection $\nabla^{{M_1}}$ on ${M_1}$. Then symmetric \textit{second fundamental form} of $F$ is given by \cite{Nore_1986}
	\begin{equation}\label{sff}
		\left( \nabla F_{\ast }\right) \left({\cal Z}_1, {\cal Z}_2\right) =\nabla_{F_\ast {\cal Z}_1}^{M_2}F_{\ast}{\cal Z}_2 - F_{\ast }\left( \nabla_{\mathcal{Z}_1}^{{M_1}} {\cal Z}_2\right)
	\end{equation}
	for ${\cal Z}_1, {\cal Z}_2 \in \Gamma (T{M_1})$, where $\nabla^{M_2}$ is the Levi-Civita connection on $M_2$. In addition, by \cite{Sahin_2010} $\left( \nabla F_{\ast }\right) \left(X, Y \right)$ is completely contained in $\left({\rm range~} F_{\ast }\right)^{\perp } {\rm for\ all}\ X, Y \in \Gamma \left( \ker F_{\ast }\right)^{\perp }$. Also, for $V \in \Gamma({\rm range ~} F_\ast)^\perp$, we have the \textit{shape operator} $\mathcal{S}$ defined as \cite[p. 188]{Sahin_book}:
	\begin{equation*}
		\nabla^{M_2}_{F_\ast X} V = -\mathcal{S}_V F_\ast X + \nabla^{F \perp}_X V.
	\end{equation*}
	Furthermore, the Gauss equation for $F$ is defined as \cite[p. 189]{Sahin_book}
	{\small \begin{eqnarray}\label{Gauss_Codazzi_RM}
			g_{2}\left( R^{{M_2}}\left( F_{\ast }Z_1, F_{\ast } Z_2\right) F_{\ast }Z_3, F_{\ast }Z_4\right) &=&g_{1}\left( R^{{M_1}}\left( Z_1, Z_2\right) Z_3, Z_4\right)+g_{2}\left( \left( \nabla F_{\ast }\right) \left( Z_1, Z_3\right), \left(\nabla F_{\ast }\right) \left(Z_2, Z_4\right) \right) \nonumber \\&&-g_{2}\left( \left( \nabla F_{\ast }\right) \left( Z_1, Z_4\right), \left(\nabla F_{\ast }\right) \left(Z_2, Z_3\right) \right),
	\end{eqnarray}}
	where $Z_i \in \Gamma \left( \ker F_{\ast}\right)^{\perp}$. Here, $R^{{M_1}}$ and $R^{{M_2}}$ denote the curvature tensors of ${M_1}$ and ${M_2}$, respectively. \\
	
	At a point $p\in {M_1}$, suppose that $\{Z_{i}\}_{i=1}^{r}$ is an orthonormal basis of the horizontal space. Then scalar curvatures $2{~\rm scal}^{{\cal H}}$ and $2{~\rm scal}^{{\cal R}}$ on the horizontal and range spaces are given, respectively, by
	\begin{equation}\label{2scalH}
		2{~\rm scal}^{{\cal H}}=\sum \limits_{i,j=1}^{r}g_{1}\left( R^{{M_1}}(Z_{i}, Z_{j})Z_{j}, Z_{i}\right), 
	\end{equation} and $$2{~\rm scal}^{{\cal R}}=\sum \limits_{i, j=1}^{r} g_{2}\left( R^{{M_2}}(F_{\ast }Z_{i}, F_{\ast }Z_{j})F_{\ast}Z_{j}, F_{\ast }Z_{i}\right).$$ Consequently, we define normalized scalar curvatures for these spaces as
	\begin{equation}\label{rhoH}
		\rho^{{\cal H}}=\frac{2{~\rm scal}^{{\cal H}}}{r\left(r-1\right) },
	\end{equation}
	and
	\begin{equation*}
		\rho^{{\cal R}}=\frac{2{~\rm scal}^{{\cal R}}}{r\left( r-1\right) }.
	\end{equation*}
	Supposing $\{V_{r+1}, \dots, V_{m_2}\}$ an orthonormal basis of $\left( {\rm range~}F_{\ast}\right)^{\perp }$ we set, 
	\begin{eqnarray*}
		B_{ij}^{{\cal H}^{\alpha }} &=&g_{2}\left( (\nabla F_{\ast})(Z_{i}, Z_{j}), V_{\alpha }\right), \quad i, j=1, \dots, r, \quad \alpha=r+1, \dots, m_2, \nonumber \\ \left\Vert B^{{\cal H}}\right\Vert^{2} &=&\sum_{i, j=1}^{r}g_{2}\left((\nabla F_{\ast})(Z_{i}, Z_{j}), (\nabla F_{\ast})(Z_{i}, Z_{j})\right), \nonumber \\ {\rm trace\, }B^{{\cal H}} &=&\sum_{i=1}^{r}(\nabla F_{\ast}) \left(Z_{i}, Z_{i}\right)\text{.}
	\end{eqnarray*}
	Then the \textit{Casorati curvature} of the horizontal space is defined as
	\begin{equation*}
		C^{{\cal H}}=\frac{1}{r}\left\Vert B^{{\cal H}}\right\Vert^{2}=\frac{1}{r} \sum_{\alpha=r+1}^{m_2} \sum_{i, j=1}^{r} \left( B_{ij}^{{\cal H}^{\alpha}}\right)^{2}.
	\end{equation*}
	
	Let $L^{{\cal H}}$ be a $k$ dimensional subspace ($k\geq 2$) of the horizontal space with orthonormal basis $\{Z_{1}, \dots, Z_{k}\}$. Then its Casorati curvature $C^{L^{\mathcal{H}}}$ is given by
	\[
	C^{L^{\mathcal{H}}}=\frac{1}{k}\sum_{\alpha=r+1}^{m_2}\sum_{i, j=1}^{k}\left( B_{ij}^{{\cal H}^{\alpha }}\right)^{2}. 
	\]
	Moreover, the \textit{normalized Casorati curvatures} $\delta_{C}^{{\cal H}}(r-1)$ and $\hat{\delta}_{C}^{{\cal H}}(r-1)$ associated with the	horizontal space at a point $p$ are given by
	\begin{equation}\label{CC_H_RM_1}
		\left[ \delta_{C}^{{\cal H}}(r-1)\right]_{p}=\frac{1}{2}C_{p}^{{\cal H}} + \frac{\left( r+1\right)}{2r} \inf \{C^{L^{\mathcal{H}}}|L^{{\cal H}}\ {\rm a\ hyperplane\ of\ }(\ker F_{\ast {p}})^\perp\}, 
	\end{equation}
	and
	\begin{equation}\label{CC_H_RM_2}
		\left[ \hat{\delta}_{C}^{{\cal H}}(r-1)\right]_{p}=2C_{p}^{{\cal H}}-\frac{\left( 2r-1\right)}{2r} \sup \{C^{L^{\mathcal{H}}}|L^{{\cal H}}\ {\rm	a\ hyperplane\ of\ }(\ker F_{\ast {p}})^\perp\}.
	\end{equation}
	
	In the sequel, for any ${\cal Z}\in \Gamma (TM_2)$, we write 
	\begin{equation}\label{decompose_GCSF_RM}
		J_{2}{\cal Z}=\phi_2 {\cal Z} = P^{\cal R}{\cal Z}+{Q}^{\cal R}{\cal Z}, 
	\end{equation}
	where ${P}^{\cal R}{\cal Z}\in \Gamma ({\rm range~}F_\ast)$, ${Q}^{\cal R}{\cal Z}\in \Gamma (
	{\rm range~} F_\ast)^{\perp }$ and ${P}^{\cal R}$ is the endomorphism on $({\rm range~}F_\ast)$ such that $\|{P}^{\cal R}\|^{2} = \sum \limits_{i,j=1}^{r}\left(g_{2}(F_{\ast }Z_{i}, J_2 F_{\ast }Z_{j})\right)^{2} = \sum \limits_{i,j=1}^{r}\left(g_{2}(F_{\ast }Z_{i}, \phi_2 F_{\ast }Z_{j})\right)^{2} = \sum \limits_{i,j=1}^{r}\left(g_{2}(F_{\ast }Z_{i}, {P}^{\cal R} F_{\ast }Z_{j})\right)^{2}$.
	
	\subsection{General Casorati inequalities for Riemannian maps}\label{subsec_2_1}
	In this subsection, we give general forms of the Casorati inequalities for the Riemannian maps as follows.
	\begin{theorem}\label{General_Inq_Thm_RM} 
		Let $F:({M_1}^{m_1}, g_{1}) \to \left( {M_2}^{m_2}, g_{2}\right)$ be a Riemannian map between Riemannian manifolds with rank $r\geq 3$. Then
		\begin{equation}\label{General_Inq_Eq_RM}
			\rho^{{\cal H}}\leq \delta_{C}^{{\cal H}}(r-1)+\rho^{{\cal R}} \quad {\rm and} \quad \rho^{{\cal H}}\leq \hat{\delta}_{C}^{{\cal H}}(r-1)+\rho^{{\cal R}}.
		\end{equation}
		The equality holds in any of the above two inequalities at a point $p\in {M_1}$ if and only if, for suitable orthonormal bases, the following hold.
		\begin{equation*}
			B_{11}^{{\cal H}^{\alpha }}=B_{22}^{{\cal H}^{\alpha}}=\cdots=B_{r-1\, r-1}^{{\cal H}^{\alpha}}=\frac{1}{2}B_{rr}^{{\cal H}^{\alpha}}, 
		\end{equation*}
		\begin{equation*}
			B_{ij}^{{\cal H}^{\alpha }}=0, \quad 1 \leq i \neq j \leq r.
		\end{equation*}
		The equality conditions can be interpreted as follows. The first condition gives $g_2 (F_\ast Z_1, \mathcal{S}_{V_{\alpha}} F_\ast Z_1) = g_2 (F_\ast Z_2, \mathcal{S}_{V_{\alpha}} F_\ast Z_2) = \cdots = g_2 (F_\ast Z_{r-1}, \mathcal{S}_{V_{\alpha}} F_\ast Z_{r-1}) = \frac{1}{2} g_2 (F_\ast Z_r, \mathcal{S}_{V_{\alpha}} F_\ast Z_r)$ with respect to all directions $(V_\alpha, \text{where}~ \alpha \in \{r+1, \dots, m_2\})$. Equivalently, there exist $(m_2-r)$ mutually orthogonal unit vector fields in $\left( {\rm range~}F_{\ast}\right)^{\perp }$ such that shape operators with respect to all directions have an eigenvalue of multiplicity $(r-1)$ and that for each $V_\alpha$ the distinguished eigendirections are the same (namely $F_\ast Z_r$). Hence, the leaves of range spaces are invariantly quasi-umbilical $\cite{DHV_2008}$. The second condition gives $g_2 (F_\ast Z_j, \mathcal{S}_{V_{\alpha}} F_\ast Z_i) = 0$ with respect to all directions $(V_\alpha, \text{where}~ \alpha \in \{r+1, \dots, m_2\}$ in  $\left( {\rm range~}F_{\ast}\right)^{\perp })$. Equivalently, the shape operator matrices become diagonal, and hence commute.
	\end{theorem}
	
	\begin{proof}
		At a point $p \in {M_1}$, let $\left\{ Z_{1}, \dots, Z_{r}\right\} $, $\{F_{\ast }Z_{1}, \dots, F_{\ast }Z_{r}\}$ and $\{V_{r+1}, \dots, V_{m_2}\}$ be orthonormal bases for $(\ker F_\ast)^\perp$, $({\rm range}~F_\ast)$ and $({\rm range}~F_\ast)^\perp$, respectively. Then by putting $Z_1=Z_4=Z_i$ and $Z_2=Z_3=Z_j$ in (\ref{Gauss_Codazzi_RM}), we obtain
		\begin{equation}\label{eq-GCI-(5)}
			2{~\rm scal}^{{\cal H}}=2{~\rm scal}^{{\cal R}}+\left\Vert {\rm trace~}B^{{\cal H}}\right\Vert^{2}-r C^{{\cal H}}.
		\end{equation}
		Now, consider the quadratic polynomial
		\begin{equation}\label{P_for_RM}
			{\cal P}^{{\cal R}}=\frac{1}{2}r(r-1){C}^{{\cal H}}+\frac{1}{2}(r^{2}-1) C^{L^{\mathcal{H}}}+2{~\rm scal}^{{\cal R}}-2{~\rm scal}^{{\cal H}}, 
		\end{equation}
		and assume that the hyperplane $L^\mathcal{H}$ is spanned by $\{Z_{1}, \dots, Z_{r-1}\}$. Then using (\ref{eq-GCI-(5)}), we obtain
		\begin{align*}
			{\cal P}^{{\cal R}}& = \frac{1}{2}r(r-1) C^{\mathcal{H}} + \frac{1}{2}(r^{2}-1) C^{L^{\mathcal{H}}} - \left\Vert {\rm trace~}B^{{\cal H}}\right\Vert^{2} + rC^{\mathcal{H}} \\& = \sum_{\alpha=r+1}^{m_2} \sum_{i=1}^{r-1} \left\{ r\left(B_{ii}^{{\cal H}^{\alpha}}\right)^{2} + (r+1) \left( B_{ir}^{{\cal H}^{\alpha}}\right)^{2}\right\} \\& + \sum_{\alpha =r+1}^{m_2}\left\{ 2\left(r+1\right) \sum_{1=i<j}^{r-1}\left(B_{ij}^{{\cal H}^{\alpha}}\right)^{2} - 2 \sum_{1=i<j}^{r} B_{ii}^{{\cal H}^{\alpha}} B_{jj}^{{\cal H}^{\alpha}} + \frac{(r-1)}{2} \left( B_{rr}^{{\cal H}^{\alpha}}\right)^{2} \right\}.
		\end{align*}
		Thus, 
		\begin{equation*}
			{\cal P}^{{\cal R}}\geq \sum_{\alpha =r+1}^{m_2}\left\{r\sum_{i=1}^{r-1}\left( B_{ii}^{{\cal H}^{\alpha }}\right)^{2}+\frac{(r-1)}{2}\left( B_{rr}^{{\cal H}^{\alpha}}\right)^{2} - 2 \sum_{1=i<j}^{r} B_{ii}^{{\cal H}^{\alpha}} B_{jj}^{{\cal H}^{\alpha }}\right\}. 
		\end{equation*}
		Let $f:{\Bbb R}^{r}\to {\Bbb R\, }$\ be a quadratic form for each $r+1 \leq \alpha \leq m_2$ such that
		\begin{equation*}
			f(B_{11}^{{\cal H}^{\alpha }}, \dots, B_{rr}^{{\cal H}^{\alpha}}) = \sum_{i=1}^{r-1}r\left( B_{ii}^{{\cal H}^{\alpha}}\right)^{2} + \frac{(r-1)}{2}\left(B_{rr}^{{\cal H}^{\alpha}}\right)^{2} - 2\sum_{1=i<j}^{r}B_{ii}^{{\cal H}^{\alpha}} B_{jj}^{{\cal H}^{\alpha }}, 
		\end{equation*}
		and constrained extremum problem $\min f$ subject to $B_{11}^{{\cal H}^{\alpha }} + \cdots + B_{rr}^{{\cal H}^{\alpha}} = k^{\alpha} \in {\Bbb R}$. Then using Lemma \ref{Lemma_Tripathi}, we get the global minimum critical point, 
		\[
		B_{11}^{{\cal H}^{\alpha }}=\cdots =B_{r-1\, r-1}^{{\cal H}^{\alpha }}=\frac{k^{\alpha }}{r+1}, \quad B_{rr}^{{\cal H}^{\alpha }}=\frac{2k^{\alpha }}{r+1}.
		\]
		Moreover, $f(B_{11}^{{\cal H}^{\alpha }}, \dots, B_{rr}^{{\cal H}^{\alpha }})=0$. Consequently, we obtain
		\begin{equation}\label{P_geq_RM}
			{\cal P}^{{\cal R}}\geq 0.
		\end{equation}
		Then (\ref{P_for_RM}) and (\ref{P_geq_RM}) yield, 
		\begin{equation}\label{eq-GCI-(9)}
			\rho^{{\cal H}}\leq \frac{1}{2}C^{\mathcal{H}}+\frac{(r+1)}{2r}C^{L^{\mathcal{H}}}+\rho^{{\cal R}}.
		\end{equation}
		In a similar manner, considering the quadratic polynomial
		\[
		{\cal Q}^{{\cal R}}=2r\left( r-1\right) C^{\mathcal{H}}-\frac{1}{2}\left(
		r-1\right) \left( 2r-1\right) C^{L^{\mathcal{H}}}+2{~\rm scal}^{
			{\cal R}}-2{~\rm scal}^{{\cal H}}, 
		\]
		we obtain
		\begin{equation}\label{Q_geq_RM}
			{\cal Q}^{{\cal R}}\geq 0, 
		\end{equation}
		and
		\begin{equation}\label{eq-GCI-(11)}
			\rho^{{\cal H}} \leq 2C^{\mathcal{H}} - \frac{(2r-1)}{2r}{C}^{L^{\cal H}} + \rho^{{\cal R}}. 
		\end{equation}
		Then (\ref{General_Inq_Eq_RM}) follows by taking the infimum in (\ref{eq-GCI-(9)}) and the supremum in (\ref{eq-GCI-(11)}) on all hyperplanes $L^{{\cal H}}$. In addition, the equality cases are true by (\ref{P_geq_RM}) and (\ref{Q_geq_RM}).
	\end{proof}
	
	\subsection{Implications for Riemannian maps to generalized complex space forms}\label{sub_sec_2_2}
	In this subsection, we obtain Casorati inequalities for Riemannian maps whose target spaces are generalized complex space forms $\left( {M_2}^{m_2}\left(c_{1}, c_{2}\right), g_{2}, J_{2}\right)$ using the general inequalities obtained in Subsection \ref{subsec_2_1}. As corollaries of this subsection, we also obtain Casorati inequalities for Riemannian maps whose target spaces are real, complex, and real K\"ahler space forms. Remarkably, we also verify these inequalities with the known Casorati inequalities for Riemannian maps whose target spaces are real and complex space forms. \\
	
	Now, we give the main result of this subsection.
	
	\begin{theorem}\label{thm_inq_gcsc_RM}
		Let $F: \left( {M_1}^{m_1}, g_{1} \right) \to \left( {M_2}^{m_2}\left(c_{1}, c_{2}\right), g_{2}, J_{2}\right)$ be a Riemannian map from a Riemannian manifold to a generalized complex space form with rank $r\geq 3$. Then 
		\begin{equation*}
			\rho^{{\cal H}} \leq \delta_{C}^{{\cal H}}(r-1) + c_{1} + \frac{3c_{2}}{r\left(r-1\right)}\Vert {P}^{\cal R}\Vert^{2} \quad {\rm and} \quad \rho^{{\cal H}}\leq \hat{\delta}_{C}^{{\cal H}}(r-1)+c_{1}+\frac{3c_{2}}{r\left( r-1\right) }\Vert {P}^{\cal R}\Vert^{2}.
		\end{equation*}
		The equality cases are the same as the equality cases of Theorem $\ref{General_Inq_Thm_RM}$.
	\end{theorem}
	
	\begin{proof}
		Since $M_2$ is a generalized complex space form, using (\ref{curvature_for_gcsf}), (\ref{hor_range_isometry}), and (\ref{decompose_GCSF_RM}) we obtain,
		\begin{eqnarray*}
			&&\sum \limits_{i,j=1}^{r} g_2(R^{M_2}(F_{\ast }Z_{i},F_{\ast }Z_{j})F_{\ast }Z_{j},F_{\ast }Z_{i}) =\sum \limits_{i,j=1}^{r} c_{1} \{ g_{1}(Z_{j},Z_{j})g_{1}(Z_{i},Z_{i})-g_{1}(Z_{i},Z_{j})g_{1}(Z_{j},Z_{i})\} \\&&+\sum \limits_{i,j=1}^{r} c_{2}\left\{ g_{2}\left( F_{\ast }Z_{i}, P^{\cal R} F_{\ast }Z_{j}\right) g_{2}\left( P^{\cal R} F_{\ast }Z_{j},F_{\ast }Z_{i}\right) -g_{2}\left( F_{\ast }Z_{j}, P^{\cal R} F_{\ast}Z_{j}\right) g_{2}\left( P^{\cal R} F_{\ast }Z_{i},F_{\ast }Z_{i}\right) \right. \\&&\left. +2g_{2}\left( F_{\ast }Z_{i}, P^{\cal R} F_{\ast }Z_{j}\right) g_{2}\left(P^{\cal R} F_{\ast }Z_{j},F_{\ast }Z_{i}\right) \right\}.
		\end{eqnarray*}
		Equivalently,
		\begin{equation*}
			2{~\rm scal}^{{\cal R}}=r\left( r-1\right) c_{1} + 3c_{2} \Vert {P}^{\cal R}\Vert^{2}.
		\end{equation*}
		Hence
		\begin{equation*}
			\rho^{\cal R}=c_{1}+\frac{3c_{2}}{r\left( r-1\right) }\Vert {P}^{\cal R}\Vert^{2}.
		\end{equation*}
		Then the required inequalities follow using the above equation in (\ref{General_Inq_Eq_RM}).
	\end{proof}
	
	We have the following straightforward corollary.
	
	\begin{corollary}\label{cor_rm_1}
		Let $F:({M_1}^{m_1}, g_{1}) \to \left({M_2}^{m_2} \left( c_{1}, c_{2} \right), g_{2}\right)$ be a Riemannian map.
		\begin{enumerate}[$1.$]
			\item If ${M_2}(c)$ is a real space form, then
			\begin{equation}\label{eq-CIRRSF}
				\rho^{{\cal H}}\leq \delta_{C}^{{\cal H}}(r-1)+c~\text{~and~}~ \rho^{{\cal H}}\leq \hat{\delta}_{C}^{{\cal H}}(r-1)+c.
			\end{equation}
			
			\item If ${M_2}(c)$ is a complex space form, then
			{\small \begin{equation}\label{eq-CIRCSF}
					\rho^{{\cal H}} \leq \delta_{C}^{{\cal H}}(r-1) + \frac{c}{4} + \frac{3c}{4r\left(r-1\right) }\Vert {P}^{\cal R}\Vert^{2}~\text{~and~}~ \rho^{{\cal H}} \leq \hat{\delta}_{C}^{{\cal H}}(r-1)+\frac{c}{4}+\frac{3c}{4r\left(r-1\right) } \Vert {P}^{\cal R}\Vert^{2}.
			\end{equation}}
			
			\item If ${M_2}(c)$ is a real K\"ahler space form, then
			{\small
				\begin{equation*}
					\rho^{{\cal H}} \leq \delta_{C}^{{\cal H}}(r-1)+\frac{\left( c+3\alpha \right)}{4} + \frac{3\left( c-\alpha \right)}{4r\left( r-1\right)} \Vert {P}^{\cal R}\Vert^{2} ~\text{~and~}~ \rho^{{\cal H}}\leq \hat{\delta}_{C}^{{\cal H}}(r-1)+\frac{\left( c+3\alpha \right)}{4}+\frac{3\left( c-\alpha \right) }{4r\left( r-1\right) }\Vert {P}^{\cal R}	\Vert^{2}.
			\end{equation*}}
		\end{enumerate}
		In addition, the equality cases are the same as the equality cases of the previous theorem.
	\end{corollary}
	
	\begin{remark}\label{rem_1}
		We observe that Eqs. $(\ref{eq-CIRRSF})$ and $(\ref{eq-CIRCSF})$ are the same as $\cite[~{\rm Eq.}~(3.10)]{LLSV}$ and $\cite[~{\rm Eq.}~(4.6)]{LLSV}$, respectively.
	\end{remark}
	
	From \cite{Sahin_2010}, we recall that for an invariant Riemannian map we have $J_2({\rm range~} F_\ast) = ({\rm range~}F_\ast)$ and for an anti-invariant Riemannian map we have $J_2({\rm range~}F_\ast) \subseteq {\rm (range~} F_\ast)^\perp$. Taking into account this information, one can get the following inequalities similar to Theorem \ref{thm_inq_gcsc_RM}.
	
	\begin{proposition}\label{thm_inv_rm_1}
		Let $F:\left( {M_1}^{m_1}, g_{1}\right) \to \left( {M_2}^{m_2}\left(c_{1}, c_{2}\right), g_{2}, J_2\right) $ be a Riemannian map from a Riemannian manifold to a generalized complex space form with rank $r\geq 3$. 
		\begin{enumerate}[$1.$]
			\item If $F$ is invariant, then
			\begin{equation*}
				\rho^{{\cal H}} \leq \delta_{C}^{{\cal H}}(r-1) + c_{1} + \frac{3c_{2}}{\left(r-1\right)} ~{\rm ~and ~}~ \rho^{{\cal H}}\leq \hat{\delta}_{C}^{{\cal H}}(r-1)+c_{1}+\frac{3c_{2}}{\left( r-1\right) }.
			\end{equation*}
			
			\item If $F$ is anti-invariant, then
			\begin{equation*}
				\rho^{{\cal H}} \leq \delta_{C}^{{\cal H}}(r-1) + c_{1} ~{\rm ~and ~}~ \rho^{{\cal H}}\leq \hat{\delta}_{C}^{{\cal H}}(r-1)+c_{1}.
			\end{equation*}
		\end{enumerate}
		The equality cases are the same as the equality cases of Theorem $\ref{General_Inq_Thm_RM}$.
	\end{proposition}
	
	\subsection{Implications for Riemannian maps to generalized Sasakian space forms}\label{sub_sec_2_3}
	In this subsection, we obtain Casorati inequalities for Riemannian maps whose target spaces are generalized Sasakian space forms $\left({M_2}^{m_2}\left(c_{1}, c_{2}, c_{3}\right), \phi_2, \xi_2, \eta_2, g_{2}\right)$ using the general inequalities obtained in Subsection \ref{subsec_2_1}. As corollaries of this subsection, we also obtain Casorati inequalities for Riemannian maps whose target spaces are Sasakian, Kenmotsu, cosymplectic, and almost $C(\alpha)$ space forms. Remarkably, we also verify these inequalities with the known Casorati inequalities for Riemannian maps whose target spaces are Sasakian and Kenmotsu space forms.\\
	
	Now, we give the main result of this subsection.
	
	\begin{theorem}\label{thm_inq_gssc_RM}
		Let $F:\left( {M_1}^{m_1}, g_{1}\right) \to \left({M_2}^{m_2}\left(c_{1}, c_{2}, c_{3}\right), \phi_2, \xi_2, \eta_2, g_{2}\right)$ be a Riemannian map from a Riemannian manifold to a generalized Sasakian space form with rank $r\geq 3$. Then
		\begin{equation*}
			\rho^{{\cal H}} \leq 
			\left\{\begin{array}{ll}
				\delta_{C}^{{\cal H}}(r-1) + c_{1} + \frac{3c_{2}}{r\left(r-1\right)} \Vert {P}^{\cal R} \Vert^{2} - \frac{2}{r} c_{3}, ~{\rm if}~ \xi_2 \in \Gamma({\rm range~}F_\ast);\\				\delta_{C}^{{\cal H}}(r-1) + c_{1} + \frac{3c_{2}}{r\left(r-1\right)} \Vert {P}^{\cal R} \Vert^{2}, ~{\rm if}~ \xi_2 \in \Gamma({\rm range~}F_\ast)^\perp, 			
			\end{array} \right.
		\end{equation*}
		and
		\begin{equation*}
			\rho^{{\cal H}} \leq 
			\left\{\begin{array}{ll}
				\hat{\delta}_{C}^{{\cal H}}(r-1) + c_{1} + \frac{3c_{2}}{r\left(r-1\right)} \Vert {P}^{\cal R} \Vert^{2} - \frac{2}{r} c_{3}, ~{\rm if}~ \xi_2 \in \Gamma({\rm range~}F_\ast);\\
				\hat{\delta}_{C}^{{\cal H}}(r-1) + c_{1} + \frac{3c_{2}}{r\left(r-1\right)} \Vert {P}^{\cal R} \Vert^{2}, ~{\rm if}~ \xi_2 \in \Gamma({\rm range~}F_\ast)^\perp.			
			\end{array} \right.
		\end{equation*}
		The equality cases are the same as the equality cases of Theorem $\ref{General_Inq_Thm_RM}$.
	\end{theorem}
	
	\begin{proof}
		Since $M_2$ is a generalized Sasakian space form, using (\ref{curvature_for_gssf}), (\ref{hor_range_isometry}), and (\ref{decompose_GCSF_RM}) we obtain,
		\begin{eqnarray*}
			&&\sum \limits_{i,j=1}^{r} g_2(R^{M_2}\left( F_{\ast }Z_{i},F_{\ast }Z_{j})F_{\ast}Z_{j},F_{\ast }Z_{i}\right)=\sum \limits_{i,j=1}^{r} c_{1}\left\{ g_{1}\left( Z_{j},Z_{j}\right)g_{1}\left( Z_{i},Z_{i}\right) -g_{1}\left( Z_{i},Z_{j}\right) g_{1}\left(Z_{j},Z_{i}\right) \right\}  \\&& + \sum \limits_{i,j=1}^{r} c_{2}\left\{ g_{2}\left( F_{\ast }Z_{i}, P^{\cal R} F_{\ast }Z_{j}\right)g_{2}\left( P^{\cal R} F_{\ast }Z_{j},F_{\ast }Z_{i}\right) -g_{2}\left( F_{\ast}Z_{j},P^{\cal R} F_{\ast }Z_{j}\right) g_{2}\left( P^{\cal R} F_{\ast }Z_{i},F_{\ast}Z_{i}\right) \right.  \\&&+\left. 2g_{2}\left( F_{\ast }Z_{i},P^{\cal R} F_{\ast }Z_{j}\right) g_{2}\left(P^{\cal R} F_{\ast }Z_{j},F_{\ast }Z_{i}\right) \right\}  \\&& + \sum \limits_{i,j=1}^{r} c_{3}\left\{ \eta_2 \left( F_{\ast }Z_{i}\right) \eta_2 \left( F_{\ast}Z_{j}\right) g_{2}\left( F_{\ast }Z_{j},F_{\ast }Z_{i}\right) -\eta_2 \left(F_{\ast }Z_{j}\right) \eta_2 \left( F_{\ast }Z_{j}\right) g_{2}\left( F_{\ast}Z_{i},F_{\ast }Z_{i}\right) \right.  \\&&\left. +g_{2}\left( F_{\ast }Z_{i},F_{\ast }Z_{j}\right) \eta_2 \left(F_{\ast }Z_{j}\right) \eta_2 \left( F_{\ast }Z_{i}\right) -g_{2}\left( F_{\ast}Z_{j},F_{\ast }Z_{j}\right) \eta_2 \left( F_{\ast }Z_{i}\right) \eta_2 \left(F_{\ast }Z_{i}\right) \right\}.
		\end{eqnarray*}
		Equivalently, 
		\begin{equation*}
			2{~\rm scal}^{{\cal R}}= 
			\left\{\begin{array}{ll}
				r\left( r-1\right) c_{1}+3c_{2}\Vert {P}^{\cal R}\Vert^{2}-2(r-1)c_{3}, ~{\rm if}~ \xi_2 \in \Gamma({\rm range~}F_\ast);\\
				r\left( r-1\right) c_{1}+3c_{2}\Vert {P}^{\cal R}\Vert^{2}, ~{\rm if}~ \xi_2 \in \Gamma({\rm range~}F_\ast)^\perp.			
			\end{array} \right.
		\end{equation*}
		Hence
		\begin{equation*}
			\rho^{\cal R}= 
			\left\{\begin{array}{ll}
				c_{1}+\frac{3c_{2}}{r\left( r-1\right)}\Vert {P}^{\cal R} \Vert^{2}-\frac{2}{r}c_{3}, ~{\rm if}~ \xi_2 \in \Gamma({\rm range~}F_\ast);\\
				c_{1}+\frac{3c_{2}}{r\left( r-1\right)}\Vert {P}^{\cal R} \Vert^{2}, ~{\rm if}~ \xi_2 \in \Gamma({\rm range~}F_\ast)^\perp.	
			\end{array} \right.
		\end{equation*}
		Then the required inequalities follow using the above equations in (\ref{General_Inq_Eq_RM}).
	\end{proof}
	
	We have the following straightforward corollary.
	
	\begin{corollary}\label{cor_rm_2}
		Let $F:({M_1}^{m_1}, g_{1})\to \left( {M_2}^{m_2}\left(c_{1}, c_{2}, c_{3}\right), \phi_2, \xi_2, \eta_2, g_{2}\right) $ be a Riemannian map.
		\begin{enumerate}[$1.$]
			\item Suppose ${M_2}(c)$ is a Sasakian space form. 
			\begin{enumerate}[$(a)$]
				\item If $\xi_2 \in \Gamma({\rm range}~F_\ast)$, then
				\begin{equation}\label{Inq_Sasakian_RM}
					\rho^{{\cal H}} \leq 
					\left\{\begin{array}{l}
						\delta_{C}^{{\cal H}}(r-1)+\frac{\left( c+3\right)}{4} +\frac{3\left( c-1\right) }{4r\left( r-1\right) }\Vert {P}^{\cal R}\Vert^{2}-\frac{\left( c-1\right) }{2r}, \\ \hat{\delta}_{C}^{{\cal H}}(r-1)+ \frac{\left( c+3\right)}{4} + \frac{3\left( c-1\right) }{4r\left( r-1\right)}\Vert {P}^{\cal R}\Vert^{2}-\frac{\left( c-1\right) }{2r}.		
					\end{array} \right.
				\end{equation}
				
				\item If $\xi_2 \in \Gamma({\rm range}~F_\ast)^\perp$, then
				\begin{equation*}
					\rho^{{\cal H}} \leq 
					\left\{\begin{array}{l}
						\delta_{C}^{{\cal H}}(r-1)+\frac{\left( c+3\right)}{4} +\frac{3\left( c-1\right) }{4r\left( r-1\right) }\Vert {P}^{\cal R}\Vert^{2}, \\				\hat{\delta}_{C}^{{\cal H}}(r-1)+ \frac{\left( c+3\right)}{4} + \frac{3\left( c-1\right) }{4r\left( r-1\right)}\Vert {P}^{\cal R}\Vert^{2}.		
					\end{array} \right.
				\end{equation*}
			\end{enumerate}
			
			\item Suppose ${M_2}(c)$ is a Kenmotsu space form. 
			\begin{enumerate}[$(a)$]
				\item If $\xi_2 \in \Gamma({\rm range}~F_\ast)$, then
				$$\rho^{{\cal H}} \leq 
				\left\{\begin{array}{l}
					\delta_{C}^{{\cal H}}(r-1)+\frac{\left( c-3\right)}{4} +\frac{3\left( c+1\right) }{4r\left( r-1\right) }\Vert {P}^{\cal R}\Vert^{2}-\frac{\left( c+1\right) }{2r}, \\ 				\hat{\delta}_{C}^{{\cal H}}(r-1)+ \frac{\left( c-3\right)}{4} + \frac{3\left( c+1\right) }{4r\left( r-1\right)}\Vert {P}^{\cal R}\Vert^{2}-\frac{\left( c+1\right) }{2r}.			
				\end{array} \right.$$
				
				\item If $\xi_2 \in \Gamma({\rm range}~F_\ast)^\perp$, then
				$$\rho^{{\cal H}} \leq 
				\left\{\begin{array}{l}
					\delta_{C}^{{\cal H}}(r-1)+\frac{\left( c-3\right)}{4} +\frac{3\left( c+1\right) }{4r\left( r-1\right) }\Vert {P}^{\cal R}\Vert^{2}, \\ \hat{\delta}_{C}^{{\cal H}}(r-1)+ \frac{\left( c-3\right)}{4} + \frac{3\left( c+1\right) }{4r\left( r-1\right)}\Vert {P}^{\cal R}\Vert^{2}.
				\end{array} \right.$$
			\end{enumerate}
			
			\item Suppose ${M_2}(c)$ is a cosymplectic space form.
			\begin{enumerate}[$(a)$]
				\item If $\xi_2 \in \Gamma({\rm range}~F_\ast)$, then
				$$\rho^{{\cal H}} \leq 
				\left\{\begin{array}{l}
					\delta_{C}^{{\cal H}}(r-1)+\frac{c}{4}+\frac{3c}{4r\left(r-1\right) }\Vert {P}^{\cal R}\Vert^{2}-\frac{c}{2r}, \\ 	\hat{\delta}_{C}^{{\cal H}}(r-1)+\frac{c}{4}+\frac{3c}{4r\left(r-1\right) }\Vert {P}^{\cal R}\Vert^{2}-\frac{c}{2r}.			
				\end{array} \right.$$
				
				\item If $\xi_2 \in \Gamma({\rm range}~F_\ast)^\perp$, then
				$$\rho^{{\cal H}} \leq 
				\left\{\begin{array}{l}
					\delta_{C}^{{\cal H}}(r-1)+\frac{c}{4}+\frac{3c}{4r\left(r-1\right) }\Vert {P}^{\cal R}\Vert^{2}, \\ 	\hat{\delta}_{C}^{{\cal H}}(r-1)+\frac{c}{4}+\frac{3c}{4r\left(r-1\right) }\Vert {P}^{\cal R}\Vert^{2}.	
				\end{array} \right.$$
			\end{enumerate}
			
			\item Suppose ${M_2}(c)$ is an almost $C(\alpha )$ space form.
			\begin{enumerate}[$(a)$]
				\item If $\xi_2 \in \Gamma({\rm range}~F_\ast)$, then
				$$\rho^{{\cal H}} \leq 
				\left\{\begin{array}{l}
					\delta_{C}^{{\cal H}}(r-1)+ \frac{\left( c+3 \alpha^{2}\right)}{4} + \frac{3\left( c-\alpha^{2}\right)}{4r\left( r-1\right) }\Vert{P}^{\cal R}\Vert^{2}-\frac{\left(c-\alpha^{2}\right) }{2r}, \\				\hat{\delta}_{C}^{{\cal H}}(r-1)+\frac{\left( c+3\alpha^{2}\right)}{4}+\frac{3\left( c-\alpha^{2}\right) }{4r\left(r-1\right) }\Vert{P}^{\cal R}\Vert^{2}-\frac{\left(c-\alpha^{2}\right) }{2r}.			
				\end{array} \right.$$
				
				\item If $\xi_2 \in \Gamma({\rm range}~F_\ast)^\perp$, then
				$$\rho^{{\cal H}} \leq 
				\left\{\begin{array}{l}
					\delta_{C}^{{\cal H}}(r-1)+ \frac{\left( c+3 \alpha^{2}\right)}{4} + \frac{3\left( c-\alpha^{2}\right)}{4r\left( r-1\right) }\Vert{P}^{\cal R}\Vert^{2}, \\				\hat{\delta}_{C}^{{\cal H}}(r-1)+\frac{\left( c+3\alpha^{2}\right)}{4}+\frac{3\left( c-\alpha^{2}\right) }{4r\left(r-1\right) }\Vert{P}^{\cal R}\Vert^{2}.			
				\end{array} \right.$$
			\end{enumerate} 
		\end{enumerate}
		In addition, the equality cases are the same as the equality cases of the previous theorem.
	\end{corollary}
	
	\begin{remark}\label{rem_3}
		In $\cite{PLS_Corrigendum}$, it is mentioned that $\xi$ in $\cite[~{\rm Theorem}~ 3.1]{PLS}$ belongs to $\Gamma({\rm range~}F_\ast)^\perp$. In addition, we observe that Eq. $(\ref{Inq_Sasakian_RM})$ is the same as $\cite[~{\rm Eq.}~(0.1)]{PLS_Corrigendum}$ and $\cite[~{\rm Eq.}~(3.1)]{PLS}$, regardless of the notation.
	\end{remark}
	
	Similarly to \cite[Definition 3.1]{GA_2021}, we recall that for an invariant Riemannian map we have $\phi_2 ({\rm range~} F_\ast) = ({\rm range~}F_\ast)$ and for an anti-invariant Riemannian map we have $\phi_2 ({\rm range~}F_\ast) \subseteq {\rm (range~} F_\ast)^\perp$. Taking into account this information, one can get the following inequalities similar to Theorem \ref{thm_inq_gssc_RM}.
	
	\begin{proposition}\label{thm_inv_rm_2}
		Let $F:\left( {M_1}^{m_1}, g_{1}\right) \to \left( {M_2}^{m_2}\left(c_{1}, c_{2}, c_{3}\right), \phi_2, \xi_2, \eta_2, g_{2}\right) $ be a Riemannian map from a Riemannian manifold to a generalized Sasakian space form with rank $r\geq 3$.
		\begin{enumerate}[$1.$]
			\item If $F$ is invariant, then
			\begin{equation*}
				\rho^{{\cal H}} \leq 
				\left\{\begin{array}{ll}
					\delta_{C}^{{\cal H}}(r-1) + c_{1} + \frac{3c_{2}}{\left(r-1\right)} - \frac{2}{r} c_{3}, ~{\rm if}~ \xi_2 \in \Gamma({\rm range~}F_\ast);\\					\delta_{C}^{{\cal H}}(r-1) + c_{1} + \frac{3c_{2}}{\left(r-1\right)}, ~{\rm if}~ \xi_2 \in \Gamma({\rm range~}F_\ast)^\perp, 
				\end{array} \right.
			\end{equation*}
			and
			\begin{equation*}
				\rho^{{\cal H}} \leq 
				\left\{\begin{array}{ll}
					\hat{\delta}_{C}^{{\cal H}}(r-1) + c_{1} + \frac{3c_{2}}{\left(r-1\right)} - \frac{2}{r} c_{3}, ~{\rm if}~ \xi_2 \in \Gamma({\rm range~}F_\ast);\\					\hat{\delta}_{C}^{{\cal H}}(r-1) + c_{1} + \frac{3c_{2}}{\left(r-1\right)}, ~{\rm if}~ \xi_2 \in \Gamma({\rm range~}F_\ast)^\perp.
				\end{array} \right.
			\end{equation*}
			
			\item If $F$ is anti-invariant, then
			\begin{equation*}
				\rho^{{\cal H}} \leq 
				\left\{\begin{array}{ll}
					\delta_{C}^{{\cal H}}(r-1) + c_{1} - \frac{2}{r} c_{3}, ~{\rm if}~ \xi_2 \in \Gamma({\rm range~}F_\ast);\\					\delta_{C}^{{\cal H}}(r-1) + c_{1}, ~{\rm if}~ \xi_2 \in \Gamma({\rm range~}F_\ast)^\perp, 
				\end{array} \right.
			\end{equation*}
			and
			\begin{equation*}
				\rho^{{\cal H}} \leq 
				\left\{\begin{array}{ll}
					\hat{\delta}_{C}^{{\cal H}}(r-1) + c_{1} - \frac{2}{r} c_{3}, ~{\rm if}~ \xi_2 \in \Gamma({\rm range~}F_\ast);\\ \hat{\delta}_{C}^{{\cal H}}(r-1) + c_{1}, ~{\rm if}~ \xi_2 \in \Gamma({\rm range~}F_\ast)^\perp.
				\end{array} \right.
			\end{equation*}
		\end{enumerate}
		The equality cases are the same as the equality cases of Theorem $\ref{General_Inq_Thm_RM}$.
	\end{proposition}
	
	\begin{remark}\label{rem_inv_map_kenmotsu}
		In particular, if $F$ is an invariant or anti-invariant Riemannian map from a Riemannian manifold to a Kenmotsu space form, then by the previous proposition, we observe that the inequalities become the same as the inequalities given in $\cite[~{\rm Corollaries} ~5.1~ {\rm and}~ 5.2]{Zaidi_Shanker}$.
	\end{remark}
	
	\section{General Casorati inequalities and implications for vertical distributions of Riemannian submersions}\label{sec_3}
	This section is divided into three subsections. The first subsection is dedicated to the general form of Casorati inequalities for the vertical distributions of Riemannian submersions between Riemannian manifolds. The second and third subsections are dedicated to the implications: Casorati inequalities when the source spaces are generalized complex and generalized Sasakian space forms, respectively. Now we start with some preliminaries.\\
	
	The Riemannian map defined in the previous section is called a \textit{Riemannian submersion} if $({\rm range}~F_\ast)^\perp = \{0\}$. Let $R^{M_1}$ and $\hat{R}$ denote the curvature tensors for $M_1$ and the vertical distribution, respectively. Then for $U_i \in \Gamma(\ker F_\ast)$, the Gauss-Codazzi equation for submersion $F$ is given as \cite{Falcitelli_2004}
	\begin{equation}\label{Gauss_Codazz_RS_T}
		g_1\left( R^{{M_1}}\left(U_1, U_2\right) U_3, U_4\right) = g_2\left( \hat{R}\left(U_1, U_2\right) U_3, U_4\right) -g_1(T_{U_1} U_4, T_{U_2}U_3) +g_1(T_{U_2} U_4, T_{U_1} U_4), 
	\end{equation}
	where $T$ is the O'Neill tensor defined in \cite{Neill_1966}.\\
	
	At a point $p \in {M_1}$, suppose that $\{U_{1}, \dots, U_{r}\}$ is an orthonormal basis of the vertical space and $\{Z_{r+1}, \dots, Z_{m_1}\}$ is an orthonormal basis of the horizontal space. Then the scalar curvatures are given by
	\begin{equation*}
		2{~\rm scal}_{M_1}^{\cal V} = \sum \limits_{i,j=1}^{r} g_{1}\left( R^{{M_1}}\left( U_{i}, U_{j}\right)U_{j}, U_{i}\right), \quad 2{~\rm scal}^{\cal V}=\sum \limits_{i,j=1}^{r} g_{1}\left(\hat{R}\left( U_{i}, U_{j}\right)U_{j}, U_{i}\right).
	\end{equation*}
	Consequently, we define normalized scalar curvatures as
	\begin{equation*}
		\rho^{\cal V}=\frac{2{~\rm scal}^{\cal V}}{r\left(r-1\right) }, \quad {\rho}_{M_1}^{\cal V}=\frac{2{~\rm scal}_{M_1}^{\cal V} }{r\left( r-1\right)}. 
	\end{equation*}
	Also, we set
	\begin{eqnarray*}
		T_{ij}^{\mathcal{V}^{\alpha}} &=& g_{1}\left(T_{U_{i}} U_{j}, Z_{\alpha}\right), \quad i, j=1, \dots, r, \quad \alpha =r+1, \dots, m_1, \nonumber \\ \left\Vert T^{\cal V}\right\Vert^{2}&=&\sum_{i, j=1}^{r}g_{1}\left(T_{U_{i}} U_{j}, T_{U_{i}} U_{j}\right), \nonumber \\ {\rm trace~}T^{\cal V} &=&\sum_{i=1}^{r}T_{U_{i}} U_{i}.
	\end{eqnarray*}
	Then the \textit{Casorati curvature} of the vertical space is defined as
	\begin{equation*}
		C^{\cal V}=\frac{1}{r}\left\Vert T^{\cal V}\right\Vert^{2}=\frac{1}{r}\sum_{\alpha=r+1}^{m_1}\sum_{i, j=1}^{r}\left( T_{ij}^{{\cal V}^{\alpha }}\right)^{2}.
	\end{equation*}
	
	Let $L^{{\cal V}}$ be a $k$ dimensional subspace ($k\geq 2$) of vertical space with orthonormal basis $\{U_{1}, \dots, U_{k}\}$. Then its Casorati curvature $C^{L^{\mathcal{V}}}$ is given by
	\begin{equation*}
		C^{L^{\cal V}} = \frac{1}{k}\sum_{\alpha =r+1}^{m_1} \sum_{i, j=1}^{k}\left(T_{ij}^{{\cal V}^{\alpha}}\right)^{2}.
	\end{equation*}
	Moreover, the \textit{normalized Casorati curvatures} $\delta_{C}^{\cal V}(r-1)$ and $\hat{\delta}_{C}^{\cal V}(r-1)$ associated with the vertical space at a point $p$ are given by
	\begin{equation*}
		\left[ \delta_{C}^{\cal V}(r-1)\right]_{p} = \frac{1}{2} C_{p}^{\cal V} + \frac{\left(r+1\right)}{2r} \inf \{C^{L^{\cal V}}|L^{\cal V}\ {\rm a\ hyperplane\ of\ } (\ker F_{\ast {p}})\}, 
	\end{equation*}
	and
	\begin{equation*}
		\left[ \hat{\delta}_{C}^{\cal V}(r-1)\right]_{p} = 2 C_{p}^{\cal V}-\frac{\left(2r-1\right)}{2r} \sup \{C^{L^{\cal V}}|L^{\cal V}\ {\rm a\ hyperplane\ of\ }(\ker F_{\ast {p}})\}.
	\end{equation*}
	
	In the sequel, for any ${\cal Z}\in \Gamma (TM_1)$, we write
	\begin{equation}\label{decompose_GCSF_VRS}
		J_1 {\cal Z}=\phi_1 {\cal Z} = P^{\cal V}{\cal Z}+{Q}^{\cal V}{\cal Z}, 
	\end{equation}
	where ${P}^{\cal V}{\cal Z}\in \Gamma ({\rm ker~}F_\ast)$, ${Q}^{\cal V}{\cal Z}\in \Gamma ({\rm ker~} F_\ast)^{\perp }$ and ${P}^{\cal V}$ is the endomorphism on $({\rm ker~}F_\ast)$ such that $\|{P}^{\cal V}\|^{2} = \sum \limits_{i,j=1}^{r}\left(g_{1}(U_{i}, J_1 U_{j})\right)^{2} = \sum \limits_{i,j=1}^{r}\left(g_{1}(U_{i}, \phi_1 U_{j})\right)^{2} = \sum \limits_{i,j=1}^{r}\left(g_{1}(U_{i}, {P}^{\cal V} U_{j})\right)^{2}$.
	
	\subsection{General Casorati inequalities for vertical distributions of Riemannian submersions}\label{subsec_3_1}
	In this subsection, we give general forms of the Casorati inequalities for the vertical distributions of Riemannian submersions as follows.
	\begin{theorem}\label{General_Inq_Thm_VD_RS}
		Let $F:({M_1}^{m_1}, g_{1}) \to \left({M_2}^{m_2}, g_{2}\right)$ be a Riemannian submersion between Riemannian manifolds with vertical space of dimension $r\geq 3$. Then
		\begin{equation}\label{eq-GCIV-(1a)}
			\rho^{\cal V} \leq \delta_{C}^{\cal V}\left( r-1\right) + {\rho}_{M_1}^{\cal V} \quad {\rm and} \quad \rho^{\cal V} \leq \hat{\delta}_{C}^{\cal V}\left( r-1\right) + {\rho}_{M_1}^{\cal V}.
		\end{equation}
		Moreover, the equality holds in any of the above two inequalities at a point $p\in {M_1}$ if and only if for suitable orthonormal bases, the following hold.
		\begin{equation*}
			T_{11}^{{\cal V}^{\alpha }}=T_{22}^{{\cal V}^{\alpha }}=\cdots =T_{r-1\, r-1}^{{\cal V}^{\alpha }}=\frac{1}{2}T_{rr}^{{\cal V}^{\alpha }}, 
		\end{equation*}
		\begin{equation*}
			T_{ij}^{{\cal V}^{\alpha }}=0, \quad 1 \leq i\neq j \leq r. 
		\end{equation*}
		The equality conditions can be interpreted as follows. The first condition gives $g_1(U_1, T_{U_1} Z_\alpha) = g_1(U_2, T_{U_2} Z_\alpha) = \cdots =g_1(U_{r-1}, T_{U_{r-1}} Z_\alpha) = \frac{1}{2}g_1(U_r, T_{U_r} Z_\alpha)$ with respect to all horizontal directions $(Z_\alpha, \text{where}~ \alpha \in \{r+1, \dots, m_1\})$. Equivalently, there exist $(m_1-r)$ mutually orthogonal horizontal unit vector fields such that the shape operator with respect to all directions has an eigenvalue of multiplicity $(r-1)$ and that for each $Z_\alpha$ the distinguished eigendirections are the same (namely $U_r$). Hence, the leaves of vertical space (called fibers of $F$) are invariantly quasi-umbilical \cite{DHV_2008}. The second condition gives $g_1(U_j, T_{U_i} Z_\alpha) = 0$ with respect to all horizontal directions $(Z_\alpha, \text{where}~ \alpha \in \{r+1, \dots, m_1\})$. Equivalently, the shape operator matrices become diagonal, and hence commute.
	\end{theorem}
	
	\begin{proof} 
		At a point $p \in {M_1}$, let $\left\{ U_{1}, \dots, U_{r}\right\}$ and $\{Z_{r+1}, \dots, Z_{m_1}\}$ be orthonormal bases for vertical and horizontal space, respectively. Then putting $U_1 = U_4 = U_i$ and $U_2 = U_3 = U_j$ in (\ref{Gauss_Codazz_RS_T}), we obtain
		\begin{equation}\label{eq-GCIV-(4)}
			2{~\rm scal}_{M_1}^{\cal V} =2{~\rm scal}^{\cal V} - \left\Vert {\rm trace~}T^{\cal V}\right\Vert^{2} + r{C}^{\cal V}.
		\end{equation}
		Now, consider the quadratic polynomial
		\begin{equation}\label{P_for_VD_RS}
			{\cal P}^{\cal V}=\frac{1}{2}r(r-1){C}^{\cal V} + \frac{1}{2}(r^{2}-1){C}^{L^{\cal V}} + 2{~\rm scal}_{M_1}^{\cal V} -2{~\rm scal}^{\cal V}, 
		\end{equation}
		and assume that the hyperplane $L^{\cal V}$ is spanned by $\{U_{1}, \dots, U_{r-1}\}$. Then using (\ref{eq-GCIV-(4)}), we obtain
		\begin{eqnarray*}
			{\cal P}^{\cal V} &=&\sum_{\alpha =r+1}^{m_1}\sum_{i=1}^{r-1}\left\{ r\left(T_{ii}^{{\cal V}^{\alpha }}\right)^{2}+(r+1)\left( T_{ir}^{{\cal V}^{\alpha }}\right)^{2}\right\} \\&& + \sum_{\alpha =r+1}^{m_1}\left\{ 2\left( r+1\right) \sum_{1=i<j}^{r-1}\left(T_{ij}^{{\cal V}^{\alpha}}\right)^{2} - 2\sum_{1=i<j}^{r}T_{ii}^{{\cal V}^{\alpha}}T_{jj}^{{\cal V}^{\alpha }}+\frac{(r-1)}{2} \left(T_{rr}^{{\cal V}^{\alpha }}\right)^{2}\right\}.
		\end{eqnarray*}
		Thus, 
		\begin{equation*}
			{\cal P}^{\cal V} \geq \sum_{\alpha =r+1}^{m_1} \left( r\sum_{i=1}^{r-1}\left(T_{ii}^{{\cal V}^{\alpha }}\right)^{2}+\frac{(r-1)}{2}\left( T_{rr}^{{\cal V}^{\alpha}}\right)^{2} - 2 \sum_{1=i<j}^{r}T_{ii}^{{\cal V}^{\alpha }}T_{jj}^{{\cal V}^{\alpha}}\right). 
		\end{equation*}
		Let $f:{\Bbb R}^{r}\to {\Bbb R\, }$\ be a quadratic form for each $r+1 \leq \alpha \leq m_1$ such that 
		\begin{equation*}
			f(T_{11}^{{\cal V}^{\alpha }}, \dots, T_{rr}^{{\cal V}^{\alpha}})=\sum_{i=1}^{r-1} r\left( T_{ii}^{{\cal V}^{\alpha }}\right)^{2} + \frac{(r-1)}{2}\left( T_{rr}^{{\cal V}^{\alpha}}\right)^{2} - 2 \sum_{1=i<j}^{r}T_{ii}^{{\cal V}^{\alpha}}T_{jj}^{{\cal V}^{\alpha }}, 
		\end{equation*}
		and constrained extremum problem $\min f$	subject to $T_{11}^{{\cal V}^{\alpha }}+\cdots +T_{rr}^{{\cal V}^{\alpha }}=k^{\alpha} \in {\Bbb R}$. Then using Lemma \ref{Lemma_Tripathi}, we get the global minimum critical point, 
		\[
		T_{11}^{{\cal V}^{\alpha }}=\cdots =T_{r-1\, r-1}^{{\cal V}^{\alpha }}=\frac{k^{\alpha }}{r+1}, \quad T_{rr}^{{\cal V}^{\alpha }}=\frac{2k^{\alpha }}{r+1}.
		\]
		Moreover, $f(T_{11}^{{\cal V}^{\alpha }}, \dots, T_{rr}^{{\cal V}^{\alpha }})=0$. Consequently, we obtain
		\begin{equation}\label{P_geq_VD_RS}
			{\cal P}^{\cal V}\geq 0.
		\end{equation}
		Then (\ref{P_for_VD_RS}) and (\ref{P_geq_VD_RS}) yield, 
		\begin{equation}\label{eq-GCIV-(10)}
			\rho^{\cal V} \leq \frac{1}{2} {C}^{\cal V} + \frac{(r+1)}{2r}{C}^{L^{\cal V}} + {\rho}_{M_1}^{\cal V}.
		\end{equation}
		In a similar manner, considering the quadratic polynomial
		\[
		{\cal Q}^{\cal V} = 2 r \left( r-1\right) {C}^{\cal V} - \frac{1}{2}\left( r-1\right) \left(2r-1\right) {C}^{L^{\cal V}} + 2{~\rm scal}_{M_1}^{\cal V} -2{~\rm scal}^{\cal V}, 
		\]
		we obtain
		\begin{equation}\label{Q_geq_VD_RS}
			{\cal Q}^{\cal V} \geq 0, 
		\end{equation}
		and
		\begin{equation}\label{eq-GCIV-(12)}
			\rho^{\cal V}\leq 2{C}^{\cal V} - \frac{\left( 2r-1\right) }{2r}{C}^{L^{\cal V}} + {\rho}_{M_1}^{\cal V}.
		\end{equation}
		Then (\ref{eq-GCIV-(1a)}) follows by taking the infimum in (\ref{eq-GCIV-(10)}) and the supremum in (\ref{eq-GCIV-(12)}) on all hyperplanes $L^{\cal V}$. In addition, the equality cases are true by (\ref{P_geq_VD_RS}) and (\ref{Q_geq_VD_RS}).
	\end{proof}
	
	\subsection{Implications for the vertical distributions of Riemannian submersions from generalized complex space forms}\label{sub_sec_3_2}
	In this subsection, we obtain Casorati inequalities for the vertical distributions of Riemannian submersions whose source spaces are generalized complex space forms $\left({M_1}^{m_1}\left( c_{1}, c_{2}\right), g_{1}, J_1\right)$ using the general inequalities obtained in Subsection \ref{subsec_3_1}. As corollaries of this subsection, we also obtain Casorati inequalities for Riemannian submersions whose source spaces are real, complex, and real K\"ahler space forms. Remarkably, we also verify these inequalities with the known Casorati inequalities for Riemannian submersions whose source spaces are real and complex space forms. \\
	
	Now, we give the main result of this subsection.
	
	\begin{theorem}\label{gcsc_rs_vert}
		Let $F: \left({M_1}^{m_1}\left( c_{1}, c_{2}\right), g_{1}, J_1\right) \to ({M_2}^{m_2}, g_{2})$ be a Riemannian submersion from a generalized complex space form onto a Riemannian manifold with vertical space of dimension $r\geq 3$. Then
		\begin{equation*}
			\rho^{\cal V}\leq \delta_{C}^{\cal V}\left( r-1\right) +c_{1}+\frac{3c_{2}}{r\left(				r-1\right) }\Vert {P}^{\cal V}\Vert^{2} \quad {\rm and} \quad \rho^{\cal V}\leq \hat{\delta}_{C}^{\cal V}\left( r-1\right) +c_{1}+\frac{3c_{2}}{r\left(				r-1\right) }\Vert {P}^{\cal V}\Vert^{2}.
		\end{equation*}
		The equality cases are the same as the equality cases of Theorem $\ref{General_Inq_Thm_VD_RS}$.
	\end{theorem}
	
	\begin{proof}
		Since $M_1$ is a generalized complex space form, using (\ref{curvature_for_gcsf}) and (\ref{decompose_GCSF_VRS}) we obtain,
		\begin{eqnarray*}
			&&\sum \limits_{i,j=1}^{r} g_1(R^{M_1}(U_{i}, U_{j})U_{j}, U_{i}) =\sum \limits_{i,j=1}^{r} c_{1} \{ g_{1}(U_{j},U_{j})g_{1}(U_{i},U_{i})-g_{1}(U_{i},U_{j})g_{1}(U_{j},U_{i})\} \\&& + \sum \limits_{i,j=1}^{r} c_{2} \left\{ g_{1}\left(U_{i}, P^{\cal V} U_{j}\right) g_{1} \left( P^{\cal V} U_{j}, U_{i}\right) - g_{1}\left( U_{j}, P^{\cal V} U_{j}\right) g_{1}\left( P^{\cal V} U_{i}, U_{i}\right) + 2 g_{1}\left( U_{i}, P^{\cal V} U_{j}\right) g_{1}\left(P^{\cal V} U_{j}, U_{i}\right) \right\}.
		\end{eqnarray*}
		Equivalently,
		\begin{equation*}
			2{~\rm scal}_{M_1}^{\cal V} = r \left( r-1\right) c_{1}+3c_{2}\Vert {P}^{\cal V}\Vert^{2}.
		\end{equation*}
		Hence
		\begin{equation*}
			{\rho}_{M_1}^{\cal V} = c_{1}+\frac{3c_{2}}{r\left( r-1\right) }\Vert {P}^{\cal V}\Vert^{2}.
		\end{equation*}
		Then the required inequalities follow using the above equation in (\ref{eq-GCIV-(1a)}).
	\end{proof}
	
	We have the following straightforward corollary.
	
	\begin{corollary}\label{cor_vrs_1}
		Let $F:\left( {M_1}^{m_1}\left( c_{1}, c_{2}\right), g_{1}\right) \to ({M_2}^{m_2}, g_{2})$ be a Riemannian submersion.
		\begin{enumerate}[$1.$]
			\item If ${M_1}(c)$ is a real space form, then
			\begin{equation}\label{eq-CIRRSF_rs}
				\rho^{{\cal V}}\leq \delta_{C}^{{\cal V}}(r-1)+c~\text{~and~}~ \rho^{{\cal V}}\leq \hat{\delta}_{C}^{{\cal V}}(r-1)+c.
			\end{equation}
			
			\item If ${M_1}(c)$ is a complex space form, then
			{\small \begin{equation}\label{eq-CIRCSF_rs}
					\rho^{{\cal V}} \leq \delta_{C}^{{\cal V}}(r-1) + \frac{c}{4} + \frac{3c}{4r\left(r-1\right) }\Vert {P}^{\cal V}\Vert^{2}~\text{~and~}~ \rho^{{\cal V}} \leq \hat{\delta}_{C}^{{\cal V}}(r-1)+\frac{c}{4}+\frac{3c}{4r\left(r-1\right) } \Vert {P}^{\cal V}\Vert^{2}.
			\end{equation}}
			
			\item If ${M_1}(c)$ is a real K\"ahler space form, then
			{\small
				\begin{equation*}
					\rho^{{\cal V}} \leq \delta_{C}^{{\cal V}}(r-1)+\frac{\left( c+3\alpha \right)}{4} + \frac{3\left( c-\alpha \right)}{4r\left( r-1\right)} \Vert {P}^{\cal V}\Vert^{2} ~\text{~and~}~ \rho^{{\cal V}}\leq \hat{\delta}_{C}^{{\cal V}}(r-1)+ \frac{\left( c+3\alpha	\right)}{4}+\frac{3\left( c-\alpha \right) }{4r\left( r-1\right) }\Vert {P}^{\cal V}	\Vert^{2}.
			\end{equation*}}
		\end{enumerate}
		In addition, the equality cases are the same as the equality cases of the previous theorem.
	\end{corollary}
	
	\begin{remark}\label{rem_2}
		We observe that Eqs. $(\ref{eq-CIRRSF_rs})$ and $(\ref{eq-CIRCSF_rs})$ are the same as $\cite[~{\rm Eq.}~(5.3)]{LLSV}$ and $\cite[~{\rm Eqs.}~(5.9) ~{\rm and}~(5.10)]{LLSV}$, respectively.
	\end{remark}
	
	From \cite{Sahin_2010a, Sahin_2013}, we recall that for an invariant Riemannian submersion we have $J_1 ({\rm ker~} F_\ast) = ({\rm ker~}F_\ast)$ and for an anti-invariant Riemannian submersion we have $J_1 ({\rm ker~}F_\ast) \subseteq {\rm (ker~} F_\ast)^\perp$. Taking into account this information, one can get the following inequalities similar to Theorem \ref{gcsc_rs_vert}.
	
	\begin{proposition}\label{thm_inv_rs_1}
		Let $F:\left( {M_1}^{m_1}\left( c_{1}, c_{2}\right), g_{1}, J_1\right) \to ({M_2}^{m_2}, g_{2})$ be a Riemannian submersion from a generalized complex space form onto a Riemannian manifold with vertical space of dimension $r\geq 3$.
		\begin{enumerate}[$1.$]
			\item If $F$ is invariant, then
			\begin{equation*}
				\rho^{{\cal V}} \leq \delta_{C}^{{\cal V}}(r-1) + c_{1} + \frac{3c_{2}}{\left(r-1\right)} ~{\rm ~and ~}~ \rho^{{\cal V}}\leq \hat{\delta}_{C}^{{\cal V}}(r-1)+c_{1}+\frac{3c_{2}}{\left( r-1\right) }.
			\end{equation*}
			
			\item If $F$ is anti-invariant, then
			\begin{equation*}
				\rho^{{\cal V}} \leq \delta_{C}^{{\cal V}}(r-1) + c_{1} ~{\rm ~and ~}~ \rho^{{\cal V}}\leq \hat{\delta}_{C}^{{\cal V}}(r-1)+c_{1}.
			\end{equation*}
		\end{enumerate}
		The equality cases are the same as the equality cases of Theorem $\ref{General_Inq_Thm_VD_RS}$.
	\end{proposition}
	
	\begin{remark}\label{rem_inv_submersion}
		In particular, if $F$ is an anti-invariant Riemannian submersion from a complex space form to a Riemannian manifold, then by the previous proposition we observe that the inequalities become the same as the inequalities given in $\cite[~{\rm Eq.~ (5.11)}]{LLSV}$.
	\end{remark}
	
	\subsection{Implications for the vertical distributions of Riemannian submersions from generalized Sasakian space forms}\label{sub_sec_3_3}
	In this subsection, we obtain Casorati inequalities for the vertical distributions of Riemannian submersions whose source spaces are generalized Sasakian space forms $\left( {M_1}^{m_1}\left( c_{1}, c_{2}, c_{3}\right), \phi_1, \xi_1, \eta_1, g_{1} \right)$ using the general inequalities obtained in Subsection \ref{subsec_3_1}. As corollaries of this subsection, we also obtain Casorati inequalities for Riemannian submersions whose source spaces are Sasakian, Kenmotsu, cosymplectic, and almost $C(\alpha)$ space forms. Remarkably, we also verify these inequalities with the known Casorati inequalities for Riemannian submersions whose source spaces are Sasakian space forms. \\
	
	Now, we give the main result of this subsection.
	
	\begin{theorem}\label{gssc_rs_vert}
		Let $F:\left( {M_1}^{m_1}\left( c_{1}, c_{2}, c_{3}\right), \phi_1, \xi_1, \eta_1, g_{1} \right) \to ({M_2}^{m_2}, g_{2})$ be a Riemannian submersion from a generalized Sasakian space form onto a Riemannian manifold with vertical space of dimension $r\geq 3$. Then
		\begin{equation*}
			\rho^{{\cal V}} \leq 
			\left\{\begin{array}{ll}
				\delta_{C}^{{\cal V}}(r-1) + c_{1} + \frac{3c_{2}}{r\left(r-1\right)} \Vert {P} ^{\cal V}\Vert^{2} - \frac{2}{r} c_{3}, ~{\rm if}~ \xi_1 \in \Gamma({\rm ker~}F_\ast);\\ \delta_{C}^{{\cal V}}(r-1) + c_{1} + \frac{3c_{2}}{r\left(r-1\right)} \Vert {P}^{\cal V} \Vert^{2}, ~{\rm if}~ \xi_1 \in \Gamma({\rm ker~}F_\ast)^\perp, 			
			\end{array} \right.
		\end{equation*}
		and
		\begin{equation*}
			\rho^{{\cal V}} \leq 
			\left\{\begin{array}{ll}
				\hat{\delta}_{C}^{{\cal V}}(r-1) + c_{1} + \frac{3c_{2}}{r\left(r-1\right)} \Vert {P}^{\cal V} \Vert^{2} - \frac{2}{r} c_{3}, ~{\rm if}~ \xi_1 \in \Gamma({\rm ker~}F_\ast);\\	\hat{\delta}_{C}^{{\cal V}}(r-1) + c_{1} + \frac{3c_{2}}{r\left(r-1\right)} \Vert {P}^{\cal V} \Vert^{2}, ~{\rm if}~ \xi_1 \in \Gamma({\rm ker~}F_\ast)^\perp.			
			\end{array} \right. 
		\end{equation*}
		The equality cases are the same as the equality cases of Theorem $\ref{General_Inq_Thm_VD_RS}$.
	\end{theorem}
	
	\begin{proof}
		Since $M_1$ is a generalized Sasakian space form, using (\ref{curvature_for_gssf}) and (\ref{decompose_GCSF_VRS}) we obtain,
		\begin{eqnarray*}
			&&\sum \limits_{i,j=1}^{r} g_1(R^{M_1}\left(U_{i}, U_{j}) U_{j}, U_{i}\right)=\sum \limits_{i,j=1}^{r} c_{1}\left\{ g_{1}\left( U_{j},U_{j}\right)g_{1}\left( U_{i},U_{i}\right) -g_{1}\left( U_{i},U_{j}\right) g_{1}\left(U_{j},U_{i}\right) \right\}  \\&& + \sum \limits_{i,j=1}^{r} c_{2}\left\{ g_{1}\left(U_{i}, P^{\cal V} U_{j}\right)g_{1}\left( P^{\cal V} U_{j}, U_{i}\right) -g_{1}\left( U_{j},P^{\cal V} U_{j}\right) g_{1}\left( P^{\cal V} U_{i},U_{i}\right) + 2g_{1}\left( U_{i},P^{\cal V} U_{j}\right) g_{1}\left(P^{\cal V} U_{j}, U_{i}\right) \right\}  \\&& + \sum \limits_{i,j=1}^{r} c_{3}\left\{ \eta_1 \left( U_{i}\right) \eta_1 \left(U_{j}\right) g_{1}\left( U_{j}, U_{i}\right) -\eta_1 \left(U_{j}\right) \eta_1 \left(U_{j}\right) g_{1}\left( U_{i},U_{i}\right) \right.  \\&&\left. +g_{1}\left( U_{i}, U_{j}\right) \eta_1 \left(U_{j}\right) \eta_1 \left(U_{i}\right) -g_{1}\left(U_{j}, U_{j}\right) \eta_1 \left(U_{i}\right) \eta_1 \left(U_{i}\right) \right\}.
		\end{eqnarray*}
		Equivalently, 
		\begin{equation*}
			2{~\rm scal}_{M_1}^{\cal V} = 
			\left\{\begin{array}{ll}
				r\left( r-1\right) c_{1}+3c_{2}\Vert {P}^{\cal V} \Vert^{2}-2(r-1)c_{3}, ~{\rm if}~ \xi_1 \in \Gamma({\rm ker~}F_\ast);\\	r\left( r-1\right) c_{1}+3c_{2}\Vert {P}^{\cal V} \Vert^{2}, ~{\rm if}~ \xi_1 \in \Gamma({\rm ker~}F_\ast)^\perp.			
			\end{array} \right.
		\end{equation*}
		Hence
		\begin{equation*}
			{\rho}_{M_1}^{\cal V}= 
			\left\{\begin{array}{ll}
				c_{1}+\frac{3c_{2}}{r\left( r-1\right)}\Vert {P}^{\cal V} \Vert^{2}-\frac{2}{r}c_{3}, ~{\rm if}~ \xi_1 \in \Gamma({\rm ker~}F_\ast);\\				c_{1}+\frac{3c_{2}}{r\left( r-1\right)}\Vert {P}^{\cal V} \Vert^{2}, ~{\rm if}~ \xi_1 \in \Gamma({\rm ker~}F_\ast)^\perp.	
			\end{array} \right.
		\end{equation*}
		Then the required inequalities follow using the above equations in (\ref{eq-GCIV-(1a)}).
	\end{proof}
	
	We have the following straightforward corollary.
	
	\begin{corollary}\label{cor_vrs_2}
		Let $F:\left( {M_1}^{m_1}\left( c_{1}, c_{2}, c_{3}\right), \phi_1, \xi_1, \eta_1, g_{1}\right) \to \left( {M_2}^{m_2}, g_{2}\right) $ be a Riemannian submersion.
		\begin{enumerate}[$1.$]
			\item Suppose ${M_1}(c)$ is a Sasakian space form. 
			\begin{enumerate}[$(a)$]
				\item If $\xi_1 \in \Gamma(\ker F_\ast)$, then
				\begin{equation*}
					\rho^{{\cal V}} \leq 
					\left\{\begin{array}{l}
						\delta_{C}^{{\cal V}}(r-1)+\frac{\left( c+3\right)}{4} +\frac{3\left( c-1\right) }{4r\left( r-1\right) }\Vert {P}^{\cal V}\Vert^{2}-\frac{\left( c-1\right) }{2r}, \\	\hat{\delta}_{C}^{{\cal V}}(r-1)+ \frac{\left( c+3\right)}{4} + \frac{3\left( c-1\right) }{4r\left( r-1\right)}\Vert {P}^{\cal V}\Vert^{2}-\frac{\left( c-1\right) }{2r}.			
					\end{array} \right.
				\end{equation*}
				
				\item If $\xi_1 \in \Gamma(\ker F_\ast)^\perp$, then
				\begin{equation}\label{Inq_Sasakian_RS}
					\rho^{{\cal V}} \leq 
					\left\{\begin{array}{l}
						\delta_{C}^{{\cal V}}(r-1)+\frac{\left( c+3\right)}{4} +\frac{3\left( c-1\right) }{4r\left( r-1\right) }\Vert {P}^{\cal V}\Vert^{2}, \\	\hat{\delta}_{C}^{{\cal V}}(r-1)+ \frac{\left( c+3\right)}{4} + \frac{3\left( c-1\right) }{4r\left( r-1\right)}\Vert {P}^{\cal V}\Vert^{2}.
					\end{array} \right.
				\end{equation}
			\end{enumerate}
			
			\item Suppose ${M_1}(c)$ is a Kenmotsu space form. 
			\begin{enumerate}[$(a)$]
				\item If $\xi_1 \in \Gamma(\ker F_\ast)$, then
				$$\rho^{{\cal V}} \leq 
				\left\{\begin{array}{l}
					\delta_{C}^{{\cal V}}(r-1)+ \frac{\left( c-3\right)}{4} +\frac{3\left( c+1\right) }{4r\left( r-1\right) }\Vert {P}^{\cal V}\Vert^{2}-\frac{\left( c+1\right) }{2r}, \\ \hat{\delta}_{C}^{{\cal V}}(r-1)+ \frac{\left( c-3\right)}{4} + \frac{3\left( c+1\right) }{4r\left( r-1\right)}\Vert {P}^{\cal V}\Vert^{2}-\frac{\left( c+1\right) }{2r}.			
				\end{array} \right.$$
				
				\item If $\xi_1 \in \Gamma(\ker F\ast)^\perp$, then
				$$\rho^{{\cal V}} \leq 
				\left\{\begin{array}{l}
					\delta_{C}^{{\cal V}}(r-1)+ \frac{\left( c-3\right)}{4} +\frac{3\left( c+1\right) }{4r\left( r-1\right) }\Vert {P}^{\cal V}\Vert^{2}, \\ \hat{\delta}_{C}^{{\cal V}}(r-1)+ \frac{\left( c-3\right)}{4} + \frac{3\left( c+1\right) }{4r\left( r-1\right)}\Vert {P}^{\cal V}\Vert^{2}.	
				\end{array} \right.$$
			\end{enumerate}

			\item Suppose ${M_1}(c)$ is a cosymplectic space form.
			\begin{enumerate}[$(a)$]
				\item If $\xi_1 \in \Gamma(\ker F_\ast)$, then
				$$\rho^{{\cal V}} \leq 
				\left\{\begin{array}{l}
					\delta_{C}^{{\cal V}}(r-1)+\frac{c}{4}+\frac{3c}{4r\left(r-1\right) }\Vert {P}^{\cal V}\Vert^{2}-\frac{c}{2r}, \\ 	\hat{\delta}_{C}^{{\cal V}}(r-1)+\frac{c}{4}+\frac{3c}{4r\left(r-1\right) }\Vert {P}^{\cal V}\Vert^{2}-\frac{c}{2r}.			
				\end{array} \right.$$
				
				\item If $\xi_1 \in \Gamma(\ker F_\ast)^\perp$, then
				$$\rho^{{\cal V}} \leq 
				\left\{\begin{array}{l}
					\delta_{C}^{{\cal V}}(r-1)+\frac{c}{4}+\frac{3c}{4r\left(r-1\right) }\Vert {P}^{\cal V}\Vert^{2}, \\ 	\hat{\delta}_{C}^{{\cal V}}(r-1)+\frac{c}{4}+\frac{3c}{4r\left(r-1\right) }\Vert {P}^{\cal V}\Vert^{2}.			
				\end{array} \right.$$
			\end{enumerate}
			
			\item Suppose ${M_1}(c)$ is an almost $C(\alpha)$ space form.
			\begin{enumerate}[$(a)$]
				\item If $\xi_1 \in \Gamma(\ker F_\ast)$, then
				$$\rho^{{\cal V}} \leq 
				\left\{\begin{array}{l}
					\delta_{C}^{{\cal V}}(r-1)+ \frac{\left( c+3 \alpha^{2}\right)}{4} + \frac{3\left( c-\alpha^{2}\right)}{4r\left( r-1\right) }\Vert{P}^{\cal V}\Vert^{2}-\frac{\left(c-\alpha^{2}\right) }{2r}, \\ \hat{\delta}_{C}^{{\cal V}}(r-1)+ \frac{\left( c+3\alpha^{2}\right)}{4} + \frac{3\left( c-\alpha^{2}\right) }{4r\left(r-1\right) }\Vert{P}^{\cal V}\Vert^{2}-\frac{\left(c-\alpha^{2}\right) }{2r}.	
				\end{array} \right.$$
				
				\item If $\xi_1 \in \Gamma(\ker F_\ast)^\perp$, then
				$$\rho^{{\cal V}} \leq 
				\left\{\begin{array}{l}
					\delta_{C}^{{\cal V}}(r-1)+ \frac{\left( c+3 \alpha^{2}\right)}{4} + \frac{3\left( c-\alpha^{2}\right)}{4r\left( r-1\right) }\Vert{P}^{\cal V}\Vert^{2}, \\ \hat{\delta}_{C}^{{\cal V}}(r-1)+ \frac{\left( c+3\alpha^{2}\right)}{4} + \frac{3\left( c-\alpha^{2}\right) }{4r\left(r-1\right) }\Vert{P}^{\cal V}\Vert^{2}.
				\end{array} \right.$$
			\end{enumerate} 
		\end{enumerate}
		In addition, the equality cases are the same as the equality cases of the previous theorem.		
	\end{corollary}
	
	\begin{remark}\label{rem_4}
		In $\cite{PLS_Corrigendum}$, it is mentioned that $\xi$ in $\cite[~{\rm Theorem}~ 4.2]{PLS}$ belongs to $\Gamma({\rm \ker~}F_\ast)^\perp$. In addition, we observe that Eq. $(\ref{Inq_Sasakian_RS})$ becomes the same as $\cite[~{\rm Eq.}~(0.3)]{PLS_Corrigendum}$ and $\cite[~{\rm Eq.}~(4.9)]{PLS}$, if we take the same notation and $(r+1)$-dimensional vertical space instead $r$.
	\end{remark}
	
	Similarly to \cite{Erken_Murathan, KP_2017}, we recall that for invariant Riemannian submersions, we have $\phi_1 ({\rm ker~} F_\ast) = ({\rm ker~}F_\ast)$ and for anti-invariant Riemannian submersions we have $\phi_1 ({\rm ker~}F_\ast) \subseteq {\rm (ker~} F_\ast)^\perp$. Taking into account this information, one can get the following inequalities similar to Theorem \ref{gssc_rs_vert}.
	
	\begin{proposition}\label{thm_inv_rs_2}
		Let $F:\left( {M_1}^{m_1}\left( c_{1}, c_{2}, c_{3}\right), \phi_1, \xi_1, \eta_1, g_{1}\right) \to ({M_2}^{m_2}, g_{2})$ be a Riemannian submersion from a generalized Sasakian space form onto a Riemannian manifold with vertical space of dimension $r\geq 3$.
		\begin{enumerate}[$1.$]
			\item If $F$ is invariant, then
			\begin{equation*}
				\rho^{{\cal V}} \leq 
				\left\{\begin{array}{ll}
					\delta_{C}^{{\cal V}}(r-1) + c_{1} + \frac{3c_{2}}{\left(r-1\right)} - \frac{2}{r} c_{3}, ~{\rm if}~ \xi_1 \in \Gamma({\rm ker~}F_\ast);\\					\delta_{C}^{{\cal V}}(r-1) + c_{1} + \frac{3c_{2}}{\left(r-1\right)}, ~{\rm if}~ \xi_1 \in \Gamma({\rm \ker~}F_\ast)^\perp, 
				\end{array} \right.
			\end{equation*}
			and
			\begin{equation*}
				\rho^{{\cal V}} \leq 
				\left\{\begin{array}{ll}
					\hat{\delta}_{C}^{{\cal V}}(r-1) + c_{1} + \frac{3c_{2}}{\left(r-1\right)} - \frac{2}{r} c_{3}, ~{\rm if}~ \xi_1 \in \Gamma({\rm ker~}F_\ast);\\					\hat{\delta}_{C}^{{\cal V}}(r-1) + c_{1} + \frac{3c_{2}}{\left(r-1\right)}, ~{\rm if}~ \xi_1 \in \Gamma({\rm ker~}F_\ast)^\perp.
				\end{array} \right.
			\end{equation*}
			
			\item If $F$ is anti-invariant, then
			\begin{equation*}
				\rho^{{\cal V}} \leq 
				\left\{\begin{array}{ll}
					\delta_{C}^{{\cal V}}(r-1) + c_{1} - \frac{2}{r} c_{3}, ~{\rm if}~ \xi_1 \in \Gamma({\rm ker~}F_\ast);\\					\delta_{C}^{{\cal V}}(r-1) + c_{1}, ~{\rm if}~ \xi_1 \in \Gamma({\rm ker~}F_\ast)^\perp, 
				\end{array} \right.
			\end{equation*}
			and
			\begin{equation*}
				\rho^{{\cal V}} \leq 
				\left\{\begin{array}{ll}
					\hat{\delta}_{C}^{{\cal V}}(r-1) + c_{1} - \frac{2}{r} c_{3}, ~{\rm if}~ \xi_1 \in \Gamma({\rm ker~}F_\ast);\\ \hat{\delta}_{C}^{{\cal V}}(r-1) + c_{1}, ~{\rm if}~ \xi_1 \in \Gamma({\rm ker~}F_\ast)^\perp.
				\end{array} \right.
			\end{equation*}
		\end{enumerate}
		The equality cases are the same as the equality cases of Theorem $\ref{General_Inq_Thm_VD_RS}$.
	\end{proposition}
	
	\section{General Casorati inequalities and implications for horizontal distributions of Riemannian submersions}\label{sec_4}
	This section is divided into three subsections. The first subsection is dedicated to the general form of Casorati inequalities for the horizontal distributions of Riemannian submersions between Riemannian manifolds. The second and third subsections are dedicated to the implications: Casorati inequalities when the source spaces are generalized complex and generalized Sasakian space forms, respectively. Now we start with some preliminaries.\\
	
	Let $F:(M_{1}^{m_{1}}, g_{1})\to \left( M_{2}^{m_{2}}, g_{2}\right) $ be a Riemannian submersion between two Riemannian manifolds. Let $R^{M_{1}}$ and $R^{\cal H}$ denote the curvature tensors for $M_{1}$ and the horizontal distribution, respectively. Then for $Z_i \in \Gamma(\ker F_\ast)^\perp$, the Gauss-Codazzi equation for submersion $F$ is given as \cite{Falcitelli_2004}
	\begin{eqnarray}\label{Gauss_Codazz_RS_A}
		R^{M_{1}}\left( Z_1, Z_2, Z_3, Z_4\right) &=&R^{\cal H}\left(Z_1, Z_2, Z_3, Z_4\right) +2 g_{1}\left( A_{Z_1} Z_2, A_{Z_3} Z_4\right)  \nonumber \\ &&-g_{1}\left( A_{Z_2} Z_3, A_{Z_1} Z_4\right) +g_{1}\left( A_{Z_1} Z_3, A_{Z_2} Z_4\right)
	\end{eqnarray}
	where $A$ is the O'Neill tensor defined in \cite{Neill_1966}.\\
	
	At a point $p\in M_{1}$, suppose that $\{Z_{i}\}_{i=1}^{r}$ is an orthonormal basis of the horizontal space and $\{U_{i}\}_{i=r+1}^{m_1}$ is an orthonormal basis of the vertical space. Then the scalar curvatures are given by (\ref{2scalH}) and $2 ~{\rm scal}_{\cal H}^{\cal H} = \sum \limits_{i,j=1}^{r} g_{1}\left( R^{\cal H}\left( Z_{i}, Z_{j}\right) Z_{j}, Z_{i}\right)$.
	Consequently, we define normalized scalar curvatures as in (\ref{rhoH}) and ${\rho}_{\cal H}^{\cal H}=\frac{2{~\rm scal}_{\cal H}^{\cal H} }{r\left( r-1\right)}$.
	Also, we set
	\begin{eqnarray*}
		A_{ij}^{\mathcal{H}^{\alpha}} &=& g_{1}\left(A_{Z_{i}} Z_{j}, U_{\alpha}\right), \quad i, j=1, \dots, r, \quad \alpha =r+1, \dots, m_1, \nonumber \\ \left\Vert A^{\cal H}\right\Vert^{2}&=&\sum_{i, j=1}^{r}g_{1}\left(A_{Z_{i}} Z_{j}, A_{Z_{i}} Z_{j}\right), \nonumber \\ {\rm trace~}A^{\cal H} &=&\sum_{i=1}^{r}A_{Z_{i}} Z_{i}.
	\end{eqnarray*}
	Then the \textit{Casorati curvature} of the horizontal space is defined as
	\begin{equation*}
		C^{\cal H}=\frac{1}{r}\left\Vert A^{\cal H}\right\Vert^{2}=\frac{1}{r}\sum_{\alpha=r+1}^{m_1}\sum_{i, j=1}^{r}\left( A_{ij}^{{\cal H}^{\alpha }}\right)^{2}.
	\end{equation*}
	
	Let $L^{{\cal H}}$ be a $k$ dimensional subspace ($k\geq 2$) of horizontal space with orthonormal basis $\{Z_{1}, \dots, Z_{k}\}$. Then its Casorati curvature $C^{L^{\mathcal{H}}}$ is given by
	\begin{equation*}
		C^{L^{\cal H}} = \frac{1}{k}\sum_{\alpha =r+1}^{m_1} \sum_{i, j=1}^{k}\left(A_{ij}^{{\cal H}^{\alpha}}\right)^{2}.
	\end{equation*}
	Moreover, the \textit{normalized Casorati curvatures} $\delta_{C}^{\cal H}(r-1)$ and $\hat{\delta}_{C}^{\cal H}(r-1)$ associated with the horizontal space at a point $p$ are given by (\ref{CC_H_RM_1}) and (\ref{CC_H_RM_2}).\\
	
	In the sequel, for any ${\cal Z}\in \Gamma (TM_1)$, we write
	\begin{equation}\label{decompose_GCSF_HRS}
		J_1 {\cal Z} = \phi_1 {\cal Z} = P^{\cal H} {\cal Z} + Q^{\cal H} {\cal Z}, 
	\end{equation}
	where $P^{\cal H} {\cal Z} \in \Gamma (\ker F_{\ast })^{\perp }$, $Q^{\cal H} {\cal Z} \in \Gamma (\ker F_{\ast})$ and $P^{\cal H}$ is the endomorphism on $(\ker F_{\ast })^{\perp}$ such that $\|{P}^{\cal H}\|^{2} = \sum \limits_{i,j=1}^{r}\left(g_{1}(Z_{i}, J_1 Z_{j})\right)^{2} = \sum \limits_{i,j=1}^{r}\left(g_{1}(Z_{i}, \phi_1 Z_{j})\right)^{2} = \sum \limits_{i,j=1}^{r}\left(g_{1}(Z_{i}, {P}^{\cal H} Z_{j})\right)^{2}$.
	
	\subsection{General Casorati inequalities for horizontal distributions of Riemannian submersions}\label{sub_sec_4_1}
	In this subsection, we give general forms of the Casorati inequalities for the horizontal distributions of Riemannian submersions as follows.
	
	\begin{theorem}\label{General_Inq_Thm_HD_RS}
		Let $F:(M_{1}^{m_{1}}, g_{1}) \to \left(M_{2}^{m_{2}}, g_{2}\right) $ be a Riemannian submersion between Riemannian manifolds with horizontal space of dimension $r\geq 3$. Then
		\begin{equation}\label{eq-GCIH-(1a)}
			\rho_{\cal H}^{\cal H} \leq \delta _{C}^{{\cal H}}\left( r-1\right) + \rho^{{\cal H}} \quad {\rm and} \quad \rho_{\cal H}^{{\cal H}}\leq \hat{\delta}_{C}^{{\cal H}} \left( r-1\right) +\rho^{{\cal H}}.
		\end{equation}
		Moreover, equality holds in any of the above two inequalities at a point $p\in M_{1}$ if and only if the tensor $A$ vanishes. Geometrically, equality means that the horizontal distribution is integrable.
	\end{theorem}
	
	\begin{proof}
		At a point $p\in M_{1}$, let $\left\{ Z_{1}, \ldots, Z_{r}\right\} $ and $\{U_{r+1}, \ldots, U_{m_{1}}\}$ be orthonormal bases for horizontal and vertical space, respectively. Then putting $Z_1=Z_4=Z_{i}$ and $Z_2=Z_3=Z_{j}$ in (\ref{Gauss_Codazz_RS_A}), we obtain
		\begin{equation}\label{eq-GCIH-(2)}
			2~{\rm scal}^{\cal H} = 2~{\rm scal}_{\cal H}^{\cal H} + 3 r{C}^{{\cal H}} - \left \Vert \text{{\rm trace}~}A^{\cal H}\right\Vert^{2}.
		\end{equation}
		Now, consider the quadratic polynomial
		\begin{equation}\label{P_for_HD_RS}
			{\cal P}^{{\cal H}}=\frac{1}{2}r(r-1){C}^{{\cal H}}+\frac{1}{2}(r^{2}-1){C}^{L^{{\cal H}}} + 2~{\rm scal}^{\cal H} - 2~{\rm scal}_{\cal H}^{\cal H}, 
		\end{equation}
		and assume that the hyperplane $L^{{\cal H}}$ is spanned by $\left\{ Z_{1}, \ldots, Z_{r-1}\right\}$. Then using (\ref{eq-GCIH-(2)}), we obtain
		\begin{eqnarray*}
			{\cal P}^{{\cal H}} &=&\sum_{\alpha=r+1}^{m_1}\left\{(r+3)\sum_{i, j=1}^{r-1} \left( A_{ii}^{{\cal H}^{\alpha }} \right)^{2} + \frac{1}{2}(r+5)\left( A_{rr}^{{\cal H}^{\alpha }}\right)^{2} +2(r+3)\sum_{1\leq i<j \leq r-1}\left( A_{ij}^{{\cal H}^{\alpha }}\right) ^{2}\right. \\&&\left. +(r+5) \sum_{i=1}^{r-1} \left( A_{ir}^{{\cal H}^{\alpha }}\right)^{2} - \sum_{i, j=1}^{r} A_{ii}^{{\cal H}^{\alpha }} A_{jj}^{{\cal H}^{\alpha }} \right\}
		\end{eqnarray*}
		Thus, 
		\begin{equation*}
			{\cal P}^{{\cal H}}\geq \sum_{\alpha = r+1}^{m_1} \left\{(r+3)\sum_{i, j=1}^{r-1}\left( A_{ii}^{{\cal H}^{\alpha }}\right)^{2} + \frac{1}{2}(r+5)\left( A_{rr}^{{\cal H}^{\alpha }}\right)^{2} -\sum_{i, j=1}^{r} A_{ii}^{{\cal H}^{\alpha }} A_{jj}^{{\cal H}^{\alpha }}\right\}.
		\end{equation*}
		Since $A_{X}Y=-A_{Y}X$, $A_{ii}^{{\cal H}^{\alpha }}=0$. Hence, 
		\begin{equation}\label{P_geq_HD_RS}
			{\cal P}^{{\cal H}}\geq 0.
		\end{equation}
		Then (\ref{P_for_HD_RS}) and (\ref{P_geq_HD_RS}) yield, 
		\begin{equation}\label{eq-GCIH-(7)}
			\rho_{\cal H}^{{\cal H}}\leq \frac{1}{2}{\cal C}^{{\cal H}}+\frac{(r+1)}{2r}{\cal C}^{L^{{\cal H}}}+\rho^{{\cal H}}.
		\end{equation}
		In a similar manner, considering the quadratic polynomial
		\begin{equation*}
			{\cal Q}^{{\cal H}}=2r\left( r-1\right) {C}^{{\cal H}}-\frac{1}{2}\left(r-1\right) \left( 2r-1\right) {C}^{L^{{\cal H}}} + 2~{\rm scal}^{\cal H} - 2~{\rm scal}_{\cal H}^{\cal H}, 
		\end{equation*}
		we obtain
		\begin{equation}\label{Q_geq_HD_RS}
			{\cal Q}^{\cal H} \geq 0, 
		\end{equation}
		and
		\begin{equation}\label{eq-GCIH-(10)}
			\rho_{\cal H}^{\cal H} \leq 2{C}^{{\cal H}} -\frac{(2r-1)}{2r}{C}^{L^{{\cal H}}} + \rho^{{\cal H}}.
		\end{equation}
		Then (\ref{eq-GCIH-(1a)}) follows by taking the infimum in (\ref{eq-GCIH-(7)}) and the supremum in (\ref{eq-GCIH-(10)}) on all hyperplanes $L^{{\cal H}}$. In addition, the equality cases are true by (\ref{P_geq_HD_RS}) and (\ref{Q_geq_HD_RS}).
	\end{proof}
	
	\subsection{Implications for the horizontal distributions of Riemannian submersions from generalized complex space forms}\label{sub_sec_4_2}
	In this subsection, we obtain Casorati inequalities for the horizontal distributions of Riemannian submersions whose source spaces are generalized complex space forms $\left( M_{1}^{m_{1}}\left( c_{1}, c_{2}\right), g_{1}, J_1\right)$ using the general inequalities obtained in Subsection \ref{sub_sec_4_1}. As corollaries of this subsection, we also obtain Casorati inequalities for Riemannian submersions whose source spaces are real, complex, and real K\"ahler space forms. \\
	
	Now, we give the main result of this subsection.
	
	\begin{theorem}\label{gcsc_rs_hor}
		Let $F:\left( M_{1}^{m_{1}}\left( c_{1}, c_{2}\right), g_{1}, J_1\right) \to (M_{2}^{m_{2}}, g_{2})$ be a Riemannian submersion from a generalized complex space form onto a Riemannian manifold with horizontal space of dimension $r\geq 3$. Then
		\begin{equation*}
			\rho_{\cal H}^{\cal H} \leq \delta_{C}^{{\cal H}}\left( r-1\right) + c_{1}+\frac{3c_{2}}{r\left( r-1\right) }\Vert {P}^{\cal H}\Vert^{2} \quad {\rm and} \quad \rho_{\cal H}^{\cal H} \leq \hat{\delta}_{C}^{{\cal H}}\left( r-1\right) + c_{1}+\frac{3c_{2}}{r\left( r-1\right) }\Vert {P}^{\cal H}\Vert^{2}.
		\end{equation*}
		The equality cases are the same as the equality cases of Theorem $\ref{General_Inq_Thm_HD_RS}$.
	\end{theorem}
	
	\begin{proof}
		Since $M_1$ is a generalized complex space form, using (\ref{curvature_for_gcsf}) and (\ref{decompose_GCSF_HRS}) we obtain,
		\begin{eqnarray*}
			&&\sum \limits_{i,j=1}^{r} g_1(R^{M_1}(Z_{i}, Z_{j})Z_{j}, Z_{i}) =\sum \limits_{i,j=1}^{r} c_{1} \{ g_{1}(Z_{j},Z_{j})g_{1}(Z_{i},Z_{i})-g_{1}(Z_{i},Z_{j})g_{1}(Z_{j},Z_{i})\} \\&& + \sum \limits_{i,j=1}^{r} c_{2} \left\{ g_{1}\left(Z_{i}, P^{\cal H} Z_{j}\right) g_{1} \left( P^{\cal H} Z_{j}, Z_{i}\right) - g_{1}\left( Z_{j}, P^{\cal H} Z_{j}\right) g_{1}\left( P^{\cal H} Z_{i}, Z_{i}\right) + 2 g_{1}\left( Z_{i}, P^{\cal H} Z_{j}\right) g_{1}\left(P^{\cal H} Z_{j}, Z_{i}\right) \right\}.
		\end{eqnarray*}
		Equivalently, 
		\begin{equation*}
			2~{\rm scal}^{\cal H} = r \left( r-1\right) c_{1} + 3 c_{2} \Vert {P}^{\cal H}\Vert^{2}.
		\end{equation*}
		Hence
		\begin{equation*}
			\rho^{{\cal H}} = c_{1} + \frac{3c_{2}}{r\left( r-1\right)} \Vert {P}^{\cal H} \Vert^{2}.
		\end{equation*}
		Then the required inequalities follow using the above equation in (\ref{eq-GCIH-(1a)}).
	\end{proof}
	
	We have the following straightforward corollary.
	
	\begin{corollary}\label{cor_hrs_1}
		Let $F:\left( {M_1}^{m_1}\left( c_{1}, c_{2}\right), g_{1}\right) \to ({M_2}^{m_2}, g_{2})$ be a Riemannian submersion.
		\begin{enumerate}[$1.$]
			\item If ${M_1}(c)$ is a real space form, then
			\begin{equation*}
				\rho_{\cal H}^{{\cal H}}\leq \delta_{C}^{{\cal H}}(r-1)+c~\text{~and~}~ \rho_{\cal H}^{{\cal H}}\leq \hat{\delta}_{C}^{{\cal H}}(r-1)+c.
			\end{equation*}
			
			\item If ${M_1}(c)$ is a complex space form, then
			\begin{equation*}
				\rho_{\cal H}^{{\cal H}} \leq \delta_{C}^{{\cal H}}(r-1) + \frac{c}{4} + \frac{3c}{4r\left(r-1\right) }\Vert {P}^{\cal H}\Vert^{2}~\text{~and~}~ \rho_{\cal H}^{{\cal H}} \leq \hat{\delta}_{C}^{{\cal H}}(r-1)+\frac{c}{4}+\frac{3c}{4r\left(r-1\right) } \Vert {P}^{\cal H}\Vert^{2}.
			\end{equation*}
			
			\item If ${M_1}(c)$ is a real K\"ahler space form, then
			{\small
				\begin{equation*}
					\rho_{\cal H}^{{\cal H}} \leq \delta_{C}^{{\cal H}}(r-1)+\frac{\left( c+3\alpha \right)}{4} + \frac{3\left( c-\alpha \right)}{4r\left( r-1\right)} \Vert {P}^{\cal H}\Vert^{2} ~\text{~and~}~ \rho_{\cal H}^{{\cal H}}\leq \hat{\delta}_{C}^{{\cal H}}(r-1)+ \frac{\left( c+3\alpha	\right)}{4}+\frac{3\left( c-\alpha \right) }{4r\left( r-1\right) }\Vert {P}^{\cal H}	\Vert^{2}.
			\end{equation*}}
		\end{enumerate}
		In addition, the equality cases are the same as the equality cases of the previous theorem.
	\end{corollary}
	
	\subsection{Implications for the horizontal distributions of Riemannian submersions from generalized Sasakian space forms}\label{sub_sec_4_3}
	In this subsection, we obtain Casorati inequalities for the horizontal distributions of Riemannian submersions whose source spaces are generalized Sasakian space forms $\left( {M_1}^{m_1}\left( c_{1}, c_{2}, c_{3}\right), \phi_1, \xi_1, \eta_1, g_{1} \right)$ using the general inequalities obtained in Subsection \ref{sub_sec_4_1}. As corollaries of this subsection, we also obtain Casorati inequalities for Riemannian submersions whose source spaces are Sasakian, Kenmotsu, cosymplectic, and almost $C(\alpha)$ space forms. \\
	
	Now, we give the main result of this subsection.
	
	\begin{theorem}\label{gssc_rs_hor}
		Let $F:\left( {M_1}^{m_1}\left( c_{1}, c_{2}, c_{3}\right), \phi_1, \xi_1, \eta_1, g_{1} \right) \to ({M_2}^{m_2}, g_{2})$ be a Riemannian submersion from a generalized Sasakian space form onto a Riemannian manifold with horizontal space of dimension $r\geq 3$. Then
		\begin{equation*}
			\rho_{\cal H}^{{\cal H}} \leq 
			\left\{\begin{array}{ll}
				\delta_{C}^{{\cal H}}(r-1) + c_{1} + \frac{3c_{2}}{r\left(r-1\right)} \Vert {P}^{\cal H} \Vert^{2} - \frac{2}{r} c_{3}, ~{\rm if}~ \xi_1 \in \Gamma({\rm ker~}F_\ast)^\perp;\\ \delta_{C}^{{\cal H}}(r-1) + c_{1} + \frac{3c_{2}}{r\left(r-1\right)} \Vert {P}^{\cal H} \Vert^{2}, ~{\rm if}~ \xi_1 \in \Gamma({\rm ker~}F_\ast), 			
			\end{array} \right.
		\end{equation*}
		and
		\begin{equation*}
			\rho_{\cal H}^{{\cal H}} \leq 
			\left\{\begin{array}{ll}
				\hat{\delta}_{C}^{{\cal H}}(r-1) + c_{1} + \frac{3c_{2}}{r\left(r-1\right)} \Vert {P}^{\cal H} \Vert^{2} - \frac{2}{r} c_{3}, ~{\rm if}~ \xi_1 \in \Gamma({\rm ker~}F_\ast)^\perp;\\	\hat{\delta}_{C}^{{\cal H}}(r-1) + c_{1} + \frac{3c_{2}}{r\left(r-1\right)} \Vert {P}^{\cal H} \Vert^{2}, ~{\rm if}~ \xi_1 \in \Gamma({\rm ker~}F_\ast).
			\end{array} \right. 
		\end{equation*}
		The equality cases are the same as the equality cases of Theorem $\ref{General_Inq_Thm_HD_RS}$.
	\end{theorem}
	
	\begin{proof}
		Since $M_1$ is a generalized Sasakian space form, using (\ref{curvature_for_gssf}) and (\ref{decompose_GCSF_HRS}) we obtain,
		\begin{eqnarray*}
			&&\sum \limits_{i,j=1}^{r} g_1(R^{M_1}\left(Z_{i}, Z_{j}) Z_{j}, Z_{i}\right)=\sum \limits_{i,j=1}^{r} c_{1}\left\{ g_{1}\left( Z_{j},Z_{j}\right)g_{1}\left( Z_{i},Z_{i}\right) -g_{1}\left( Z_{i},Z_{j}\right) g_{1}\left(Z_{j},Z_{i}\right) \right\}  \\&& + \sum \limits_{i,j=1}^{r} c_{2}\left\{ g_{1}\left(Z_{i}, P^{\cal H} Z_{j}\right)g_{1}\left( P^{\cal H} Z_{j}, Z_{i}\right) -g_{1}\left( Z_{j},P^{\cal H} Z_{j}\right) g_{1}\left( P^{\cal H} Z_{i},Z_{i}\right) + 2g_{1}\left( Z_{i},P^{\cal H} Z_{j}\right) g_{1}\left(P^{\cal H} Z_{j}, Z_{i}\right) \right\}  \\&& + \sum \limits_{i,j=1}^{r} c_{3}\left\{ \eta_1 \left( Z_{i}\right) \eta_1 \left(Z_{j}\right) g_{1}\left( Z_{j}, Z_{i}\right) -\eta_1 \left(Z_{j}\right) \eta_1 \left(Z_{j}\right) g_{1}\left( Z_{i},Z_{i}\right) \right.  \\&&\left. +g_{1}\left( Z_{i}, Z_{j}\right) \eta_1 \left(Z_{j}\right) \eta_1 \left(Z_{i}\right) -g_{1}\left(Z_{j}, Z_{j}\right) \eta_1 \left(Z_{i}\right) \eta_1 \left(Z_{i}\right) \right\}.
		\end{eqnarray*}
		Equivalently, 
		\begin{equation*}
			2{~\rm scal}^{\cal H} = 
			\left\{\begin{array}{ll}
				r\left( r-1\right) c_{1}+3c_{2}\Vert {P}^{\cal H}\Vert^{2}-2(r-1)c_{3}, ~{\rm if}~ \xi_1 \in \Gamma({\rm ker~}F_\ast)^\perp;\\	r\left( r-1\right) c_{1}+3c_{2}\Vert {P}^{\cal H}\Vert^{2}, ~{\rm if}~ \xi_1 \in \Gamma({\rm ker~}F_\ast).			
			\end{array} \right.
		\end{equation*}
		Hence
		\begin{equation*}
			{\rho}^{\cal H} = 
			\left\{\begin{array}{ll}
				c_{1}+\frac{3c_{2}}{r\left( r-1\right)}\Vert {P}^{\cal H} \Vert^{2}-\frac{2}{r}c_{3}, ~{\rm if}~ \xi_1 \in \Gamma({\rm ker~}F_\ast)^\perp;\\				c_{1}+\frac{3c_{2}}{r\left( r-1\right)}\Vert {P}^{\cal H} \Vert^{2}, ~{\rm if}~ \xi_1 \in \Gamma({\rm ker~}F_\ast).	
			\end{array} \right.
		\end{equation*}
		Then the required inequalities follow using the above equations in (\ref{eq-GCIH-(1a)}).
	\end{proof}
	
	We have the following straightforward corollary.
	
	\begin{corollary}\label{cor_hrs_2}
		Let $F:\left( {M_1}^{m_1}\left( c_{1}, c_{2}, c_{3}\right), \phi_1, \xi_1, \eta_1, g_{1}\right) \to \left( {M_2}^{m_2}, g_{2}\right) $ be a Riemannian submersion.
		\begin{enumerate}[$1.$]
			\item Suppose ${M_1}(c)$ is a Sasakian space form. 
			\begin{enumerate}[$(a)$]
				\item If $\xi_1 \in \Gamma(\ker F_\ast)^\perp$, then
				$$\rho_{\cal H}^{{\cal H}} \leq 
				\left\{\begin{array}{l}
					\delta_{C}^{{\cal H}}(r-1)+\frac{\left( c+3\right)}{4} +\frac{3\left( c-1\right) }{4r\left( r-1\right) }\Vert {P}^{\cal H}\Vert^{2}-\frac{\left( c-1\right) }{2r}, \\	\hat{\delta}_{C}^{{\cal H}}(r-1)+ \frac{\left( c+3\right)}{4} + \frac{3\left( c-1\right) }{4r\left( r-1\right)}\Vert {P}^{\cal H}\Vert^{2}-\frac{\left( c-1\right) }{2r}.
				\end{array} \right.$$
				
				\item If $\xi_1 \in \Gamma(\ker F_\ast)$, then
				$$\rho_{\cal H}^{{\cal H}} \leq 
				\left\{\begin{array}{l}
					\delta_{C}^{{\cal H}}(r-1)+\frac{\left( c+3\right)}{4} +\frac{3\left( c-1\right) }{4r\left( r-1\right) }\Vert {P}^{\cal H}\Vert^{2}, \\	\hat{\delta}_{C}^{{\cal H}}(r-1)+ \frac{\left( c+3\right)}{4} + \frac{3\left( c-1\right) }{4r\left( r-1\right)}\Vert {P}^{\cal H}\Vert^{2}.
				\end{array} \right.$$
			\end{enumerate}
			
			\item Suppose ${M_1}(c)$ is a Kenmotsu space form. 
			\begin{enumerate}[$(a)$]
				\item If $\xi_1 \in \Gamma(\ker F_\ast)^\perp$, then
				$$\rho_{\cal H}^{{\cal H}} \leq 
				\left\{\begin{array}{l}
					\delta_{C}^{{\cal H}}(r-1)+ \frac{\left( c-3\right)}{4} +\frac{3\left( c+1\right) }{4r\left( r-1\right) }\Vert {P}^{\cal H}\Vert^{2}-\frac{\left( c+1\right) }{2r}, \\ \hat{\delta}_{C}^{{\cal H}}(r-1)+ \frac{\left( c-3\right)}{4} + \frac{3\left( c+1\right) }{4r\left( r-1\right)}\Vert {P}^{\cal H}\Vert^{2}-\frac{\left( c+1\right) }{2r}.
				\end{array} \right.$$
				
				\item If $\xi_1 \in \Gamma(\ker F_\ast)$, then
				$$\rho_{\cal H}^{{\cal H}} \leq 
				\left\{\begin{array}{l}
					\delta_{C}^{{\cal H}}(r-1)+ \frac{\left( c-3\right)}{4} +\frac{3\left( c+1\right) }{4r\left( r-1\right) }\Vert {P}^{\cal H}\Vert^{2}, \\ \hat{\delta}_{C}^{{\cal H}}(r-1)+ \frac{\left( c-3\right)}{4} + \frac{3\left( c+1\right) }{4r\left( r-1\right)}\Vert {P}^{\cal H}\Vert^{2}.			
				\end{array} \right.$$
			\end{enumerate}
			
			\item Suppose ${M_1}(c)$ is a cosymplectic space form.
			\begin{enumerate}[$(a)$]
				\item If $\xi_1 \in \Gamma(\ker F_\ast)^\perp$, then
				$$\rho_{\cal H}^{{\cal H}} \leq 
				\left\{\begin{array}{l}
					\delta_{C}^{{\cal H}}(r-1)+\frac{c}{4}+\frac{3c}{4r\left(r-1\right) }\Vert {P}^{\cal H}\Vert^{2}-\frac{c}{2r}, \\ 	\hat{\delta}_{C}^{{\cal H}}(r-1)+\frac{c}{4}+\frac{3c}{4r\left(r-1\right) }\Vert {P}^{\cal H}\Vert^{2}-\frac{c}{2r}.			
				\end{array} \right.$$
				
				\item If $\xi_1 \in \Gamma(\ker F_\ast)$, then
				$$\rho_{\cal H}^{{\cal H}} \leq 
				\left\{\begin{array}{l}
					\delta_{C}^{{\cal H}}(r-1)+\frac{c}{4}+\frac{3c}{4r\left(r-1\right) }\Vert {P}^{\cal H}\Vert^{2}, \\ 	\hat{\delta}_{C}^{{\cal H}}(r-1)+\frac{c}{4}+\frac{3c}{4r\left(r-1\right) }\Vert {P}^{\cal H}\Vert^{2}.			
				\end{array} \right.$$
			\end{enumerate}
			
			\item Suppose ${M_1}(c)$ is an almost $C(\alpha)$ space form.
			\begin{enumerate}[$(a)$]
				\item If $\xi_1 \in \Gamma(\ker F_\ast)^\perp$, then
				$$\rho_{\cal H}^{{\cal H}} \leq 
				\left\{\begin{array}{l}
					\delta_{C}^{{\cal H}}(r-1)+ \frac{\left( c+3 \alpha^{2}\right)}{4} + \frac{3\left( c-\alpha^{2}\right)}{4r\left( r-1\right) }\Vert{P}^{\cal H}\Vert^{2}-\frac{\left(c-\alpha^{2}\right) }{2r}, \\ \hat{\delta}_{C}^{{\cal H}}(r-1)+ \frac{\left( c+3\alpha^{2}\right)}{4} + \frac{3\left( c-\alpha^{2}\right) }{4r\left(r-1\right) }\Vert{P}^{\cal H}\Vert^{2}-\frac{\left(c-\alpha^{2}\right) }{2r}.
				\end{array} \right.$$
				
				\item If $\xi_1 \in \Gamma(\ker F_\ast)$, then
				$$\rho_{\cal H}^{{\cal H}} \leq 
				\left\{\begin{array}{l}
					\delta_{C}^{{\cal H}}(r-1)+ \frac{\left( c+3 \alpha^{2}\right)}{4} + \frac{3\left( c-\alpha^{2}\right)}{4r\left( r-1\right) }\Vert{P}^{\cal H}\Vert^{2}, \\ \hat{\delta}_{C}^{{\cal H}}(r-1)+ \frac{\left( c+3\alpha^{2}\right)}{4} + \frac{3\left( c-\alpha^{2}\right) }{4r\left(r-1\right) }\Vert{P}^{\cal H}\Vert^{2}.	
				\end{array} \right.$$
			\end{enumerate} 
		\end{enumerate}
		In addition, the equality cases are the same as the equality cases of the previous theorem.		
	\end{corollary}
	
	\section{Examples}
	In this section, we provide explicit examples illustrating the general inequalities. \\
	
	We start with the following example of a Riemannian map for which the inequalities derived in Theorem \ref{General_Inq_Thm_RM} hold strictly.
	\begin{example}
		Let $\left( M_1 =\mathbb{R}^4, ~g_1 = \sum \limits_{i=1}^{4} (dx_i)^2 +x_1 (dx_1) (dx_2)\right)$ and $\left( M_2 =\mathbb{R}^4, ~g_2 = \sum \limits_{j=1}^{4} (dy_j)^2 + y_1 (dy_1) (dy_2)\right)$ be two Riemannian manifolds with $x_1, y_1 \in (-1, 1)$. Define a map $F: (M_1, g_1) \to (M_2, g_2)$ such that 
		$$F(x_1, x_2, x_3, x_4) = \left(x_1, 0, x_3, x_4\right).$$ Then we get
		\begin{equation*}
			(\ker F_\ast) = \operatorname{span}\left\{U_1 = \frac {\partial}{\partial x_2}\right\},
		\end{equation*}
		\begin{equation*}
			(\ker F_\ast)^\perp = \operatorname{span}\left\{Z_1= \frac {\partial}{\partial x_1}, Z_2 =\frac {\partial}{\partial x_3} , Z_3 = \frac {\partial}{\partial x_4}\right\},
		\end{equation*}
		\begin{equation*}
			({\rm range}~ F_\ast) = \operatorname{span}\left\{F_\ast (Z_1) = \frac{\partial}{\partial y_1}, F_\ast (Z_2) = \frac{\partial}{\partial y_3}, F_\ast (Z_3) = \frac{\partial}{\partial y_4}\right\},
		\end{equation*}
		and
		\begin{equation*}
			({\rm range}~ F_\ast)^\perp = \operatorname{span}\left\{V_1 = \frac{\partial}{\partial y_2} \right\},
		\end{equation*}
		where $\left\{\frac{\partial}{\partial x_i} \right\}_{i=1}^{4}$ and $\left\{\frac{\partial}{\partial y_j} \right\}_{j=1}^{4}$ are bases of $T_p M_1$ and $T_{F(p)}M_2$, respectively. For all $i, j \in \{1,2,3\}$, we verify that $$g_1 (Z_i, Z_j) = g_2 (F_\ast Z_i, F_\ast Z_j).$$ Therefore, $F$ is a Riemannian map with rank $3$. The only non-zero Christoffel symbols for $g_1$ are $\Gamma^{x_1}_{x_1,x_1} = \frac{x_1}{x_1^2 -1}, \Gamma^{x_2}_{x_1,x_1} = \frac{-1}{x_1^2 -1}$, and for $g_2$ are $\Gamma^{y_1}_{y_1,y_1} = \frac{y_1}{y_1^2 -1}, \Gamma^{y_2}_{y_1,y_1} = \frac{-1}{y_1^2 -1}$. Then by some computations using all Christoffel symbols, Eq. $(\ref{sff})$ and the fact that $y_1=x_1$, we obtain
		\begin{equation*}
			(\nabla F_\ast) (Z_1, Z_1) = \frac{-1}{y_1^2 -1} \frac{\partial}{\partial y_2}, (\nabla F_\ast) (Z_2, Z_2) = 0, (\nabla F_\ast) (Z_3, Z_3) =0,
		\end{equation*}
		\begin{equation*}
			(\nabla F_\ast) (Z_i, Z_j) = 0 ~ \text{for all $i \neq j \in \{1,2,3\}$}.
		\end{equation*}
		Hence, for $V_\alpha =V_1$, we have
		\begin{equation*}
			B_{11}^{{\cal H}^{\alpha }}=g_{2}\left( (\nabla F_{\ast})(Z_{1}, Z_{1}), V_{\alpha }\right)=\frac{-1}{y_1^2 -1}.
		\end{equation*}
		Moreover, in this case, we have
		\begin{equation*}
			B_{11}^{{\cal H}^{\alpha }}\neq B_{22}^{{\cal H}^{\alpha}}= 0=B_{33}^{{\cal H}^{\alpha}}, 
		\end{equation*}
		\begin{equation*}
			B_{ij}^{{\cal H}^{\alpha }}=0, \quad 1 \leq i \neq j \leq 3.
		\end{equation*}
		Clearly, for $F$ the inequalities obtained in Theorem $\ref{General_Inq_Thm_RM}$ do not achieve equality.
	\end{example}
	
	\vspace{0.2cm}
	
	Now, we give another example of a Riemannian map for which the inequalities derived in Theorem \ref{General_Inq_Thm_RM} hold with equality.
	
	\begin{example}
		Define $F: \left( M_1 =\mathbb{R}^4, g_1 = \sum \limits_{i=1}^{4} (dx_i)^2\right) \to \left( M_2 =\mathbb{R}^4, g_2 = \sum \limits_{j=1}^{4} (dy_j)^2\right)$ such that $$F(x_1, x_2, x_3, x_4) = \left(\sqrt{x_1^2 + x_2^2}, x_3, x_4, 0\right).$$ Then we get
		\begin{equation*}
			(\ker F_\ast) = \operatorname{span}\left\{U_1 = \frac {x_2} { \sqrt{x_1^2 + x_2^2}} \frac {\partial}{\partial x_1}-\frac {x_1} { \sqrt{x_1^2 + x_2^2}} \frac {\partial}{\partial x_2}\right\},
		\end{equation*}
		\begin{equation*}
			(\ker F_\ast)^\perp = \operatorname{span}\left\{Z_1= \frac {x_1} {\sqrt{x_1^2 + x_2^2}} \frac {\partial}{\partial x_1}+\frac {x_2} {\sqrt{x_1^2 + x_2^2}} \frac {\partial}{\partial x_2} ,  Z_2= \frac {\partial}{\partial x_3}, Z_3 = \frac {\partial}{\partial x_4}\right\},
		\end{equation*}
		\begin{equation*}
			({\rm range}~ F_\ast) = \operatorname{span}\left\{F_\ast (Z_j) = \frac{\partial}{\partial y_j} \right\}_{j=1}^{3}, ~ \text{and ~} ({\rm range}~ F_\ast)^\perp = \operatorname{span}\left\{V_1 = \frac{\partial}{\partial y_4} \right\},
		\end{equation*}
		where $\left\{\frac{\partial}{\partial x_i} \right\}_{i=1}^{4}$ and $\left\{\frac{\partial}{\partial y_j} \right\}_{j=1}^{4}$ are bases of $T_p M_1$ and $T_{F(p)}M_2$, respectively. One can verify that $F$ is a Riemannian map between Riemannian manifolds with rank $3$. Further, similar to the previous example, 
		for $V_\alpha =V_1$, we have $$B_{ij}^{{\cal H}^{\alpha }}=g_{2}\left( (\nabla F_{\ast})(Z_{i}, Z_{j}), V_{\alpha }\right)=0, \quad 1 \leq i, j \leq 3.$$
		Clearly, for $F$ the inequalities obtained in Theorem $\ref{General_Inq_Thm_RM}$ achieve equality.
	\end{example}
	
	\vspace{0.2cm}
	
	Now, we give an example of a Riemannian submersion for which the inequalities derived in Theorem \ref{General_Inq_Thm_VD_RS} and Theorem \ref{General_Inq_Thm_HD_RS} hold strictly and with equality, respectively.
	
	\begin{example}  
		Let $(M_1 = \{(x_1,x_2,x_3,x_4,x_5, x_6) \in \Bbb{R}^{6}:x_1\neq 0, x_2\neq 0, x_3\neq 0, x_4\neq 0, x_5\neq 0, x_6 \neq 0\}, g_1 = \sum \limits_{i=1}^{6} dx_i^2 )$ and $(M_2 = \{(y_1, y_2, y_3) \in \Bbb{R}^3 : y_1 >0, y_2 >0, y_3 >0\}, g_2 = \sum \limits_{j=1}^{3} dy_j^2)$ be two Riemannian manifolds. Define a map $F : (M_1, g_1) \rightarrow (M_2, g_2) $ such that
		\begin{equation*}
			F(x_1,x_2,x_3,x_4,x_5,x_6)= \left( \sqrt{x_1^2 + x_2^2}, \sqrt{x_3^2 + x_4^2} , \sqrt{x_5^2 + x_6^2} \right).
		\end{equation*}
		Fixing $\tau_1 =  \sqrt{x_1^2 + x_2^2}$, $ \tau_2 =  \sqrt{x_3^2 + x_4^2}$, and $ \tau_3 =  \sqrt{x_5^2 + x_6^2}$, we get
		\begin{equation*}
			(\ker F_\ast) = \operatorname{span}
			\left\{ U_1= \frac {x_2} {\tau_1} \frac {\partial}{\partial x_1}-\frac {x_1} {\tau_1} \frac {\partial}{\partial x_2} ,  U_2= \frac {x_4} {\tau_2} \frac {\partial}{\partial x_3}-\frac {x_3} {\tau_2} \frac {\partial}{\partial x_4}, U_3 = \frac {x_6} {\tau_3} \frac {\partial}{\partial x_5}-\frac {x_5} {\tau_3} \frac {\partial}{\partial x_6} \right\},
		\end{equation*}
		\begin{equation*}
			(\ker F_\ast)^\bot = \operatorname{span} \left\{ Z_1= \frac {x_1} {\tau_1} \frac {\partial}{\partial x_1}+\frac {x_2} {\tau_1} \frac {\partial}{\partial x_2} ,  Z_2= \frac {x_3} {\tau_2} \frac {\partial}{\partial x_3}+\frac {x_4} {\tau_2} \frac {\partial}{\partial x_4}, Z_3 = \frac {x_5} {\tau_3} \frac {\partial}{\partial x_5}+\frac {x_6} {\tau_3} \frac {\partial}{\partial x_6} \right\},
		\end{equation*}
		and
		\begin{equation*}
			({\rm range}~ F_\ast) = \operatorname{span}\left\{F_\ast (Z_j) = \frac{\partial}{\partial y_j} \right\}_{j=1}^{3},
		\end{equation*}
		where $\left\{\frac{\partial}{\partial x_i}\right\}_{i=1}^{6}$ and $\left\{\frac{\partial}{\partial y_j}\right\}_{j=1}^{3}$ are bases of $T_p M_1$ and $T_{F (p)} M_2$ respectively. One can verify that $F$ is a Riemannian submersion with horizontal and vertical spaces of dimensions $3$. Further, by some computations using all zero Christoffel symbols, we obtain
		\begin{equation*}
			T_{U_1} U_1= - \left( \frac {x_1} {\tau_1^2} \frac {\partial}{\partial x_1} + \frac {x_2} {\tau_1^2} \frac {\partial}{\partial x_2} \right),
		\end{equation*}
		\begin{equation*}
			T_{U_2} U_2 = - \left( \frac {x_3} {\tau_2^2} \frac {\partial}{\partial x_3} + \frac {x_4} {\tau_2^2} \frac {\partial}{\partial x_4} \right),
		\end{equation*}
		\begin{equation*}
			T_{U_3} U_3 = - \left( \frac {x_5} {\tau_3^2} \frac {\partial}{\partial x_5} + \frac {x_6} {\tau_3^2} \frac {\partial}{\partial x_6} \right),
		\end{equation*}
		and
		\begin{equation*}
			T_{U_i} U_j= 0 ~\text{ for all $i \neq j$}.
		\end{equation*}
		Moreover,
		\begin{equation*}
			g_1(T_{U_1} U_1, Z_1) = - \frac{1}{\tau_1}, g_1(T_{U_2} U_2, Z_1) = g_1(T_{U_3} U_3, Z_1) = 0,
		\end{equation*}
		\begin{equation*}
			g_1(T_{U_1} U_1, Z_2) = 0 , g_1(T_{U_2} U_2, Z_2) = - \frac{1}{\tau_2}, g_1(T_{U_3} U_3, Z_2) = 0,
		\end{equation*}
		and
		\begin{equation*}
			g_1(T_{U_1} U_1, Z_3) = 0, g_1(T_{U_2} U_2, Z_3) = 0, g_1(T_{U_3} U_3, Z_3) = - \frac{1}{\tau_3}.
		\end{equation*}
		Equivalently,
		\begin{equation*}
			\text{for $Z_\alpha = Z_1$, we have ~}T_{11}^{\mathcal{V}^\alpha} \neq T_{22}^{\mathcal{V}^\alpha} = T_{33}^{\mathcal{V}^\alpha},
		\end{equation*}
		\begin{equation*}
			\text{for $Z_\alpha = Z_2$, we have ~}T_{11}^{\mathcal{V}^\alpha} \neq T_{22}^{\mathcal{V}^\alpha} \neq  T_{33}^{\mathcal{V}^\alpha},
		\end{equation*}
		and
		\begin{equation*}
			\text{for $Z_\alpha = Z_3$, we have ~} T_{11}^{\mathcal{V}^\alpha} = T_{22}^{\mathcal{V}^\alpha} \neq T_{33}^{\mathcal{V}^\alpha}.
		\end{equation*}
		Clearly, for $F$ the inequalities obtained in Theorem $\ref{General_Inq_Thm_VD_RS}$ do not attain equality.\\
		Furthermore, suppose that $a_i, b_i$ are some real numbers for $1 \leq i \leq 3$. Then for horizontal vectors $X=a_1 Z_1 + a_2 Z_2 + a_3 Z_3$ and $Y = b_1 Z_1 + b_2 Z_2 + b_3 Z_3$, we obtain
		\begin{equation*}
			A_X Y = \frac{1}{2} \mathcal{V}[X, Y]= \frac{1}{2}\mathcal{V}\left( \nabla^{M_1}_X Y - \nabla^{M_1}_Y X\right)= 0.
		\end{equation*}
		Therefore, $A$ vanishes; equivalently, the horizontal distribution is integrable.	Clearly, for $F$ the inequalities obtained in Theorem $\ref{General_Inq_Thm_HD_RS}$ achieve equality.
	\end{example}
	
	\vspace{0.2cm}
	In the following example, we see that the inequalities derived in Theorems $\ref{General_Inq_Thm_VD_RS}$, $\ref{General_Inq_Thm_HD_RS}$, $\ref{gssc_rs_vert}$, and $\ref{gssc_rs_hor}$ hold with equality.

	\begin{example}
		Let $f$ be a positive smooth function on $\mathbb{R}$ and let $\mathbb{C}$ be a generalized complex space form of real dimension $6$. Considering a warped product manifold $\left(M_1 = \mathbb{R} \times_f \mathbb{C}, ~g_1 = dt^2 + f^2 \sum \limits_{i=1}^{3} ((dx_i)^2 + (dy_i)^2)\right)$ and a Riemannian manifold $\left( M_2 =\mathbb{R}^4, g_2 = f^2 \sum \limits_{j=1}^{4} (dz_j)^2\right)$, we define $F: (M_1, g_1) \to (M_2, g_2)$ such that $$F(t, x_1, y_1, x_2, y_2, x_3, y_3) = (x_2, y_2, x_3, y_3).$$ Then we get
		\begin{equation*}
			(\ker F_\ast) = \operatorname{span}\left\{U_1 = \frac{\partial}{\partial t}, U_2 = \frac{1}{f} \frac{\partial}{\partial x_1}, U_3 = \frac{1}{f} \frac{\partial}{\partial y_1}\right\},
		\end{equation*}
		\begin{equation*}
			(\ker F_\ast)^\perp = \operatorname{span}\left\{Z_1 = \frac{1}{f} \frac{\partial}{\partial x_2}, Z_2 = \frac{1}{f} \frac{\partial}{\partial y_2}, Z_3 = \frac{1}{f} \frac{\partial}{\partial x_3}, Z_4 = \frac{1}{f} \frac{\partial}{\partial y_3}\right\},
		\end{equation*}
		and
		\begin{equation*}
			({\rm range}~ F_\ast) = \operatorname{span}\left\{F_\ast (Z_j) = \frac{1}{f} \frac{\partial}{\partial z_j} \right\}_{j=1}^{4},
		\end{equation*}
		where $\left\{\frac{\partial}{\partial t}, \frac{1}{f} \frac{\partial}{\partial x_1}, \frac{1}{f} \frac{\partial}{\partial y_1}, \frac{1}{f} \frac{\partial}{\partial x_2}, \frac{1}{f} \frac{\partial}{\partial y_2}, \frac{1}{f} \frac{\partial}{\partial x_3}, \frac{1}{f} \frac{\partial}{\partial y_3} \right\}$ and $\left\{\frac{\partial}{\partial z_j} \right\}_{j=1}^{4}$ are bases of $T_p M_1$ and $T_{F(p)}M_2$, respectively. One can verify that $F$ is a Riemannian submersion between Riemannian manifolds with a horizontal space of dimension $4$ and a vertical space of dimension $3$. By $\cite{ABC}$, we also see that $M_1$ is a generalized Sasakian space form with $c_1 = - \frac{{f'}^2}{f^2}$, $c_2 = 0$, and $c_3 =- \frac{{f'}^2}{f^2} + \frac{f''}{f}$. Further, we compute non-zero Christoffel symbols for the metric $g_1$, which are $\Gamma^t_{i, i} = -f f'$ and $\Gamma^i_{t, i} = \frac{f'}{f}$ for all $i \in \{x_j, y_j\}_{j=1}^{3}$. Using these values, we obtain
		\begin{equation*}
			\nabla^{M_1}_{U_1} U_i = \nabla^{M_1}_{U_2} U_3 = \nabla^{M_1}_{U_3} U_2 = 0 ~ \text{for all $1 \leq i \leq 3$},
		\end{equation*}
		\begin{equation*}
			\nabla^{M_1}_{U_2} U_1 = \frac{f'}{f^2} \frac{\partial}{\partial x_1}, \nabla^{M_1}_{U_3} U_1 = \frac{f'}{f^2} \frac{\partial}{\partial y_1}, \nabla^{M_1}_{U_2} U_2 = \nabla^{M_1}_{U_3} U_3 = -\frac{f'}{f} \frac{\partial}{\partial t},
		\end{equation*}
		\begin{equation*}
			\nabla^{M_1}_{Z_i} Z_j = 0 ~ \text{for all $1  \leq i \neq j \leq 4$ and ~} \nabla^{M_1}_{Z_i} Z_i = -\frac{f'}{f} \frac{\partial}{\partial t} ~ \text{for all $1  \leq i \leq 4$}.
		\end{equation*}
		This implies that $\mathcal{H} \nabla^{M_1}_{U_i} U_j = T_{U_i} U_j = 0$ for all $1 \leq i, j \leq 3$. Hence, for all $Z_\alpha \in \{Z_1, Z_2, Z_3, Z_4\}$, we have $$T_{ij}^{{\cal V}^{\alpha }}=g_1(T_{U_i} U_j, Z_{\alpha}) = 0, \quad 1 \leq i, j \leq 3.$$
		Clearly, for $F$ the inequalities obtained in Theorems $\ref{General_Inq_Thm_VD_RS}$ and $\ref{gssc_rs_vert}$ achieve equality. In addition, similar to the previous example, one can verify that $A$ vanishes; equivalently, the horizontal distribution is integrable. Clearly, for $F$ the inequalities obtained in Theorems $\ref{General_Inq_Thm_HD_RS}$ and $\ref{gssc_rs_hor}$ achieve equality.
	\end{example}
	
	\vspace{0.2cm}
	Lastly, we have the following example for which the inequalities derived in Theorem \ref{General_Inq_Thm_HD_RS} hold strictly. 
	
	\begin{example}
		The submersion from $\mathbb{S}^{2n+1}$ onto $\mathbb{C} P^n$ defined in $\cite[{\rm ~Section ~5}]{Neill_1966}$ has a non-integrable horizontal distribution of dimension $2n$. Clearly, for this submersion, the inequalities obtained in Theorem $\ref{General_Inq_Thm_HD_RS}$ do not attain equality when $n \geq 2$.
	\end{example}
	
	\begin{remark}
		By following similar steps to the previous examples, one can also construct Riemannian maps (resp. submersions) for which the inequalities obtained 
		in Theorems $\ref{thm_inq_gcsc_RM}$, $\ref{thm_inq_gssc_RM}$, $\ref{gcsc_rs_vert}$, and $\ref{gcsc_rs_hor}$ hold either strictly or with equality. For the existence of a Riemannian map (resp. submersion) with target (resp. source) space as a particular generalized complex space form, we refer to $\cite[{\rm ~Example ~2}]{LLSV}$ (resp. $\cite[{\rm ~Example~ 3}]{LLSV}$).
	\end{remark}
	
	\section{Concluding remarks}
	In this section, we give some remarks to the readers.\\
	
	The following remark explains the geometric meaning of Casorati invariants and the significance of derived inequalities.
	\begin{remark}
		Casorati invariants provide a geometric meaning related to the visualization of shapes and appearances. Precisely, in the theories of Riemannian maps and Riemannian submersions, these invariants indicate deviation from flatness in the orthogonal complementary spaces, such as fibers, leaves of horizontal distributions, and leaves of range spaces. For a surely curved surface, the Gaussian curvature may vanish, but the Casorati curvature will not vanish. Thus, geometrically, the Casorati curvature provides a better intuition than the Gaussian curvature. In the previous sections, we have generalized some optimal inequalities involving normalized Casorati curvatures or Casorati invariants for Riemannian maps and Riemannian submersions, and have shown that these inequalities yield relationships between normalized scalar curvatures and normalized Casorati curvatures. The main results are established by an optimization technique using curvature expressions, tensor decompositions, an algebraic result (Lemma $\ref{Lemma_Tripathi}$), etc., and showing a parabolic quadratic polynomial in the components of the second fundamental form or O'Neill tensors. Our generalizations for Riemannian maps and Riemannian submersions are valuable to geometric optimizations, enabling Casorati optimization across various manifolds. In addition, these provide relationships that address the comparison problem between intrinsic and extrinsic invariants, as the Casorati and scalar invariants are considered extrinsic and intrinsic, respectively. Moreover, their equality cases characterize specific geometrically significant submanifolds (i.e., leaves of distributions), such as integrable, invariantly quasi-umbilical, etc.
	\end{remark}
	
	The following remark describes why the structural condition $r \geq 3$ on the dimensions of vertical and horizontal distributions has been assumed.
	
	\begin{remark}
		We note that the inequalities obtained in the previous sections involve normalized Casorati curvatures for the hyperplanes of vertical and horizontal distributions. Hence, to establish the corresponding meaningful inequalities, it is mandatory to take dimensions of these distributions as $r \geq 2$. In addition, by $\cite[\text{Remark}~ 2.1]{Zhang_Zhang}$, the normalized Casorati curvatures vanish trivially for $r=2$. It is known that due to isometry, the rank of a Riemannian map is the same as the dimension of its horizontal distribution. Therefore, avoiding trivial circumstances, we have assumed $r \geq 3$.
	\end{remark}
	
	The following remark discusses the importance of generalized complex and generalized Sasakian structures in the present study.
	
	\begin{remark}
		In the previous sections, we have applied our main results to the settings of generalized complex and generalized Sasakian structures. Since these settings unify various structures under two umbrellas, they facilitate the study of common properties and derivations across these structures. For example, we have applied the derivations by direct substitutions to obtain special cases in Corollaries $\ref{cor_rm_1}$, $\ref{cor_rm_2}$, $\ref{cor_vrs_1}$, $\ref{cor_vrs_2}$, $\ref{cor_hrs_1}$, and $\ref{cor_hrs_2}$. These settings also prepare the reader to understand a theory that synthesizes existing and more specific cases. In addition, these settings make it easier for researchers familiar with the fundamentals to transition to more advanced/generalized theories. Moreover, such generalized manifolds have several applications in geometry, general relativity, theoretical physics, and related areas, including the construction of new interesting examples. Hence, the obtained relations for generalized complex and generalized Sasakian structures may serve as bridges between those areas and may attract readers to future research in this direction.
	\end{remark}

	\noindent R. Singh\\
	Department of Mathematics, Banaras Hindu University, Varanasi, Uttar Pradesh-221005, India.\\ E-mail: khandelrs@bhu.ac.in;	ORCID: 0009-0009-1270-3831\\
	
	\noindent K. Meena\\
	Department of Mathematics, Indian Institute of Technology Jodhpur, Rajasthan-342030, India.\\ E-mail: kirankapishmeena@gmail.com; ORCID: 0000-0002-6959-5853\\
	
	\noindent K. C. Meena\\
	Scientific Analysis Group, Defence Research and Development Organisation, Delhi-110054, India.\\
	E-mail: meenakapishchand@gmail.com; ORCID: 0000-0003-0182-8822\\
\end{document}